%% file: bgfc.tex
\documentclass[12pt, reqno, a4paper]{amsart}
\usepackage{graphics,epsfig}
\usepackage{amsfonts,amsmath,amssymb,amsbsy,amsthm}
\usepackage{bm}
\usepackage{color}
\usepackage{float}
\usepackage{mathrsfs}
\usepackage{bbm}
\usepackage[normalem]{ulem}

\usepackage[left=3cm,right=3cm,top=3cm,bottom=3cm]{geometry}

\usepackage{environments}

\input{notation}

\def\Oma{\Omega^{\rm a}}
\def\Omb{\Omega^{\rm b}}
\def\Omc{\Omega^{\rm c}}

\definecolor{cocol}{rgb}{0.8, 0, 0}
\definecolor{lzcol}{rgb}{0,0,0.8}

\begin{document}

\title[Atomistic/Continuum Blending with Ghost Force Correction]{
  Atomistic/Continuum Blending  \\
  with Ghost Force Correction}

\author[C. Ortner]{Christoph Ortner}
\address{Christoph Ortner\\ Mathematics Institute \\ Zeeman Building \\
  University of Warwick \\ Coventry CV4 7AL \\ UK}
\email{christoph.ortner@warwick.ac.uk}

\author{L. Zhang}
\address{L. Zhang \\ Department of Mathematics,
  Institute of Natural Sciences, and MOE Key Lab in Scientific and Engineering Computing\\
  Shanghai Jiao Tong University \\ 800 Dongchuan Road \\ Shanghai
  200240 \\ China}
\email{lzhang@sjtu.edu.cn}

\numberwithin{equation}{section}
\numberwithin{theorem}{section}

\date{\today}

\thanks{CO's work was supported by EPSRC grant EP/H003096, ERC
  Starting Grant 335120 and by the Leverhulme Trust through a Philip
  Leverhulme Prize. LZ's work was supported by One Thousand Plan of
  China for Young Scientists. }



\begin{abstract}
  We combine the ideas of atomistic/continuum energy blending and
  ghost force correction to obtain an energy-based atomistic/continuum
  coupling scheme which has, for a range of benchmark problems, the
  same convergence rates as optimal force-based coupling schemes.

  We present the construction of this new scheme, numerical results
  exploring its accuracy in comparison with established schemes, as
  well as a rigorous error analysis for an instructive special case.
  %
  % The basic idea is to modify the interatomic potential in a way that
  % does not change the energy, but significantly reduces the effect of
  % the ``ghost forces''. This construction is related to and sheds new
  % light on the popular ghost-force removal iteration proposed by {\it
  %   Shenoy et al} (1999), but remains conservative and admits a
  % rigorous analysis.
\end{abstract}

\maketitle

%===========================================================================

\section{Introduction}
\label{sec:intro}
Atomistic/continuum coupling schemes are a popular class of
multi-scale methods for concurrent coupling between atomistic and
continuum mechanics in the simulation of crystalline solids. An
ongoing effort to develop a rigorous numerical analysis, summarised in
\cite{2013-atc.acta}, has led to a distillation of many key ideas in
the field, and a number of improvements.

It remains an open problem to construct a general and practical
QNL-type a/c coupling scheme along the lines of
\cite{Shimokawa:2004,E:2006,OrZh:2013}. The present paper develops an
a/c coupling approach combining two popular {\em practical} schemes,
ghost force correction \cite{Shenoy:1999a} and blending
\cite{XiBe:2004}, as well as its rigorous analysis.

In the ghost force correction method of Shenoy et al
\cite{Shenoy:1999a} the spurious interface forces caused by a/c
coupling (usually termed {\em ghost forces}) are removed by adding a
suitable dead load correction, which can be computed either from a
previous step in a quasi-static process or through a self-consistent
iteration. In the latter case, this process was shown to be formally
equivalent to the force-based quasicontinuum (QCF) scheme
\cite{Dobson:2008a}. Difficult open problems remain in the analysis of
the QCF method; see \cite{Lu.bqcf:2011, Dobson:arXiv0903.0610,
  LuMing12, MakrOrtSul:qcf.nonlin} for recent advances, however in
practice the scheme seems to be optimal (in terms of the error
committed by the coupling mechanism). The main drawback is that the
QCF forces are non-conservative, i.e., there is no associated energy
functional.

An alternative scheme to reduce the effect of ghost forces is the
blending method of Xiao and Belytschko \cite{XiBe:2004}. Instead of a
sharp a/c interface, the atomistic and continuum models are blended
smoothly. This does not remove, but reduce the error due to ghost
forces \cite{VkLusk2011:blended, LiOrShVK:2014}. The blending variant
that we will consider is the BQCE scheme, formulated in
\cite{VkLusk2011:blended, 2012-CMAME-optbqce} and analysed in
\cite{VkLusk2011:blended, LiOrShVK:2014}. These analyses demonstrate
that altough the error due to ghost forces is reduced, it still
remains the dominant error contribution, i.e., the {\em
  bottleneck}. Moreover, if one were to generalize the BQCE scheme to
multi-lattices, then the reduced symmetries imply that the scheme
would not be convergent (in the sense of \cite{LiOrShVK:2014}).

In the present paper, we show how a variant of the ghost force
correction idea can be used to improve on the BQCE scheme and to
result in an a/c coupling that we denote blended ghost force
correction (BGFC), which is quasi-optimal within the context of the
framework developed in \cite{LiOrShVK:2014}.  Importantly, within the
context in which we present this scheme, our formulation does {\em
  not} yield the force-based BQCF scheme \cite{Lu.bqcf:2011,
  LiOrShVK:2014, 2013-CMAME-bqcfcomp} but is {\em energy-based}. As a
matter of fact, it is most instructive to construct the scheme not
from the point of view of ghost-force correction, but through a
modification of the site energies.

The remainder of the article is structured as follows. In
\S~\ref{sec:scheme} we formulate the BGFC scheme for point defects. In
\S~\ref{sec:numerical_tests} we extend the benchmark tests from
\cite{2012-CMAME-optbqce, 2013-CMAME-bqcfcomp, OrZh:2013} to the new scheme. In
\S~\ref{sec:scheme:ext} we briefly describe extensions to higher-order
finite element, dislocations and multi-lattices. Finally, in
\S~\ref{sec:analysis} we present a rigorous error analysis for the
simple-lattice point defect case.

\def\BGFC{BGFC~}

\section{The \BGFC Scheme}
\label{sec:scheme}
For the sake of simplicity of presentation we first present the \BGFC
scheme in an infinite lattice setting and for point defects. This also
allows us to focus on the benchmark problems discussed in
\cite{2013-CMAME-bqcfcomp, 2012-CMAME-optbqce, OrZh:2013} which are the main
motivation for the introduction of our new scheme. We present 
extensions to other problems in \S~\ref{sec:scheme:ext}.

\subsection{Atomistic model}
\label{sec:scheme:atm}
\def\Rdef{R^{\rm def}}
\def\Rg{\mathcal{R}}
\def\Rgnn{\mathcal{N}}
\def\rcut{r_{\rm cut}}
Consider a {\em homogeneous reference lattice} $\Lhom := \mA \Z^d$, $d
\in \{2, 3\}$, $\mA \in \R^{d \times d}$, non-singular. Let $\L
\subset \R^d$ be the {\em reference configuration}, satisfying $\L
\setminus B_{\Rdef} = \Lhom \setminus B_{\Rdef}$ and $\# (\L \cap
B_{\Rdef}) < \infty$, for some radius $\Rdef > 0$. Here, and
  throughout, $B_R := \{ x \in \R^d \sep |x| \leq R\}$.

The mismatch between $\L$ and $\Lhom$ in $B_{\Rdef}$ represents a
possible defect. For example, $\L = \Lhom$ for an impurity, $\L
  \subset \Lhom$ for a vacancy and $\L \supset \Lhom$ for an
  interstitial.

For $a \in \L$, let $\Rgnn_a \subset \L \setminus \{a\}$ be a set of
``nearest-neighbour'' directions satisfying ${\rm span}(\Rgnn_a) =
\R^d$ and $\sup_{a \in \L} \# \Rgnn_a < \infty$.  

A {\em deformed configuration} is a map $y : \L \to \R^d$. If it is
clear from the context what we mean, then we denote $r_{ab} := |y(a) -
y(b)|$ for $a, b \in \L$. We denote the identity map by $x$.

For each $a \in \L$ let $\Phi_a(y)$ denote the {\em site energy}
associated with the lattice site $a \in \L$. For example, in the EAM
model \cite{Daw:1984a},
\begin{equation}
  \label{eq:eam_potential}
  \Phi_a(y) := \sum_{b \in \L\setminus\{a\}} \phi(r_{ab}) + G\B(
  {\textstyle \sum_{b \in \L \setminus \{a\}}} \rho(r_{ab}) \B),
\end{equation}
for a pair potential $\phi$, an electron density function $\rho$ and
an embedding function $G$. We assume that the potentials have a finite
interaction range, that is, there exists $\rcut > 0$ such that
$\pp_{y(b)}^j \Phi_a(y) = 0$ for all $j \geq 1$, whenever $r_{ab} >
\rcut$, and that they are homogeneous outside the defect region
$B_{\Rdef}$. The latter assumption is discussed in detail in
\S~\ref{sec:scheme:hom}.

To describe, e.g., impurity defects, we allow $\phi, G, \rho$ to be
species-dependent, i.e., $\phi = \phi_{ab}, G = G_a, \psi =
\psi_{ab}$.

Under suitable conditions on the site potentials $\Phi_a, a \in \L$
(most crucially, regularity and homogeneity outside the defect;
cf. \S~\ref{sec:scheme:hom}), it is shown in
\cite{EhrOrtSha:defectsV1} that the energy-difference functional
\begin{equation}
  \label{eq:scheme:Ea}
  \Ea(u) := \sum_{a \in \L} \Phi'_a(u), \qquad \text{where } \Phi'_a(u)
  := \Phi_a(x+u)-\Phi_a(x),
\end{equation}
is well-defined for all relative {\em displacements} $u \in \Us$,
where $\Us$ is given by
\begin{align*}
  \Us &:= \b\{ v : \L \to \R^d \bsep |v|_{\Us} < +\infty \b\}, \quad
  \text{where} \\
  |v|_{\Us} &:= \B(\sum_{a \in \L} \sum_{b \in \Rgnn_a} |v(b)-v(a)|^2 \B)^{1/2}.
\end{align*}

\def\ua{u^\a}

The atomistic problem is to compute
\begin{equation}
  \label{eq:atm}
  \ua \in \arg\min \b\{ \Ea(v) \bsep v \in \Us \b\}.
\end{equation}
This problem is analysed in considerable detail in
\cite{EhrOrtSha:defectsV1}.

\subsection{Homogeneous crystals}
\label{sec:scheme:hom}
\def\Ldef{\L_{\rm def}}
\def\Ldefhom{\L^{\rm hom}_{\rm def}}
\def\Lren{L^{\rm ren}} 
We briefly discuss homogeneity of site potentials, a concept which is
important for the introduction of the \BGFC scheme, and also plays a
crucial role in the definition of the atomistic model \eqref{eq:atm}.

Homogeneity of the site potential outside the defect core entails
simply that only one atomic species occurs. For finite-range
interactions, this can be formalised by requiring that, for $|a|, |b|$
sufficiently large, $\Phi_a(y) = \Phi_b(z)$ where $z(\ell) =
y(\ell+a-b)$ within the interaction range of $b$, and suitably
extended outside.

For the case $\L = \Lhom$ we say that the site potentials are {\em
  globally homogeneous} if $\Phi_a(y) = \Phi_b(z)$ for all $a, b \in
\Lhom$. In this case, it is easy to see (summation by parts, or
point-symmetry of the lattice) that
\begin{displaymath}
  \sum_{a \in \Lhom} \b\< \del \Phi_a(x), u \b\> = 0 \qquad \forall
  u : \Lhom \to \R^d, {\rm supp}(u) \text{~compact}.
\end{displaymath}
  
For general $\L$, suppose that the site potentials are homogeneous
outside $B_{\Rdef}$, with finite interaction range. Let $\Phi_a^{\rm
  hom}, a \in \Lhom$ be a globally homogeneous site potential so that
$\< \del \Phi_a(x), u \> = \< \del \Phi_a^{\rm hom}, u \>$ for all $a
\in \L, |a| > \Rdef+\rcut$ and for all $u : \L \cup \Lhom \to \R^d$
with compact support.
% \lz{LZ: shall we remark that $\L = \Lhom$,
%   impurity, $\L\supset \Lhom$ interstitial, $\L\subset \Lhom$ vacancy,
%   etc?}. 
Further, for $u : \L \to \R^d$ let $Eu$ denote an arbitrary
extension to $\Lhom$ (e.g., $Eu(\ell) = 0$ for $\ell \in \Lhom
\setminus \L$), and let $\Ldef := \L \cap B_{\Rdef + \rcut}$, then
\begin{align} 
  \notag
  \Ea(u) &= \sum_{a \in \L} \Phi_a'(u) - \sum_{a \in \Lhom} \<
  \del\Phi_a^{\rm hom}(x), Eu \> \\
  \label{eq:scheme:Ea_renormalised}
  &= \sum_{a \in \L} \Phi_a''(u) + \< \Lren, Eu \>,
\end{align}
where
\begin{align}
  \label{eq:scheme:defn_Phi''}
  \Phi_a''(u) &:= \Phi_a(x+u) - \Phi_a(x) - \< \del\Phi_a(x), Eu \>, \\
  \label{eq:scheme:defn_Lren}
  \< \Lren, u \> &:= \sum_{a \in \Ldef} \<
  \del\Phi_a(x), u \> - \sum_{a \in \Ldefhom} \<
  \delta\Phi_a^{\rm hom}(x), Eu \>.
\end{align}

% \lz{LZ: missing $\delta$ in front of $\<\Phi_a^{\rm hom}(x), Eu \>$} 
The ``renormalisation'' \eqref{eq:scheme:Ea_renormalised} of $\Ea$ is
the basis for proving that $\Ea$ is well-defined on the energy space
$\Us$ \cite{EhrOrtSha:defectsV1}.

\subsection{Continuum model}
\label{sec:scheme:cb}
Suppose that the site potentials are homogeneous outside the defect
core and let $\Phi_a^{\rm hom}$ 
% \lz{LZ: do you mean $\Phi^{\rm hom}_0$ in the following formula?}
% CO: it does not matter at which site $a$ it is evaluated
 be the associated globally homogeneous
site potentials for the homogeneous lattice;
cf. \S~\ref{sec:scheme:hom}.

To formulate atomistic to continuum coupling schemes we require a
continuum model compatible with \eqref{eq:scheme:Ea} defined through a
strain energy function $W : \R^{d \times d} \to \R$. A typical choice
in the multi-scale context is the Cauchy--Born model, which is defined
via
\begin{displaymath}
  W(\mF) := \det \mA^{-1} \Phi^{\rm hom}_0(\mF x),
\end{displaymath}
which represents the energy per unit volume in the homogeneous crystal
$\mF \Lhom = \mF \mA \Z^d$. The associated strain energy difference is
denoted by $W'(\mG) := W(\mI+\mG) - W(\mI)$.

\subsection{Standard blending scheme}
\label{sec:scheme:bqce}
\def\Th{\mathcal{T}_h}
\def\Nh{\mathcal{X}_h}
\def\Ush{\Us_h}
\def\Ra{R^\a}
\def\Rb{R^{\rm b}}
\def\Rc{R^\c}
\def\Eb{\E^{\rm b}}
\def\dof{{\rm DOF}}
\begin{figure}
  \begin{center}
    \includegraphics[width=\textwidth]{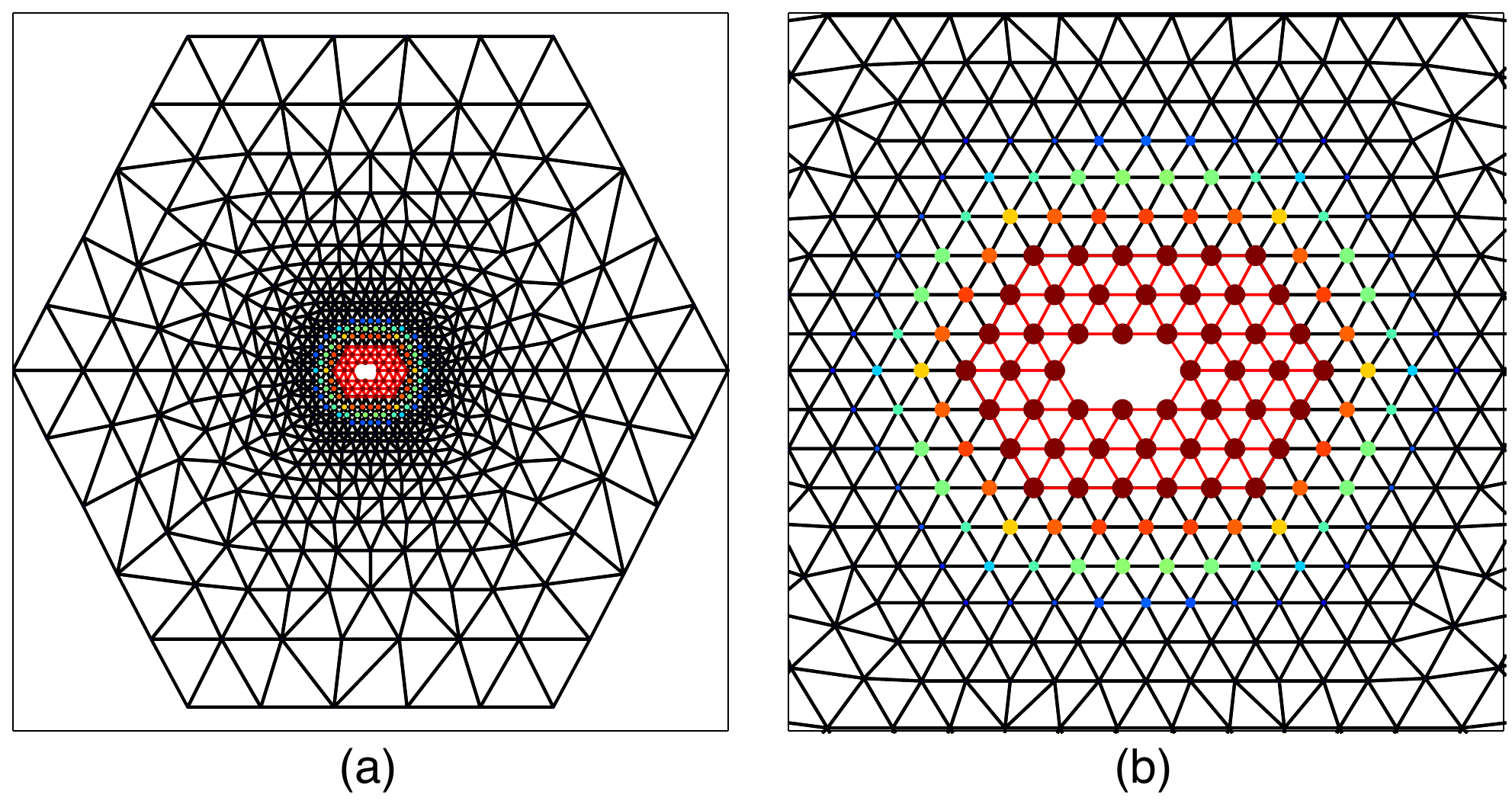}
    \caption{\label{geometry} Computational domain, finite element
      grid and atomistic region as used in the construction of the
      BQCE and BGFC schemes. The size and color of the spheres in (b)
      indicate the value of the blending function (large/red stands
      for $\beta = 0$).}
  \end{center}
\end{figure}
To formulate the blended quasicontinuum (BQCE) scheme as introduced in
\cite{VkLusk2011:blended,2012-CMAME-optbqce} and analysed in
\cite{LiOrShVK:2014} we begin by defining a regular simplicial finite
element grid $\Th$ with nodes $\Nh$, with the minimal requirement that
$\Nh \cap B_{\Rdef} = \L \cap B_{\Rdef}$ (that is, the defect core is
resolved exactly). Let $\dof := \#\Nh$.  Let $\Omh := \bigcup \Th
\supset B_{\Rc}$, with $\Rc \geq \Rdef$ be the resulting computational
domain and let the space of coarse-grained admissible displacements be
given by
\begin{displaymath}
  \Ush := \b\{ v_h \in C(\R^d; \R^d) \bsep v_h \text{ is p.w. affine
    w.r.t. } \Th, \text{ and } v_h|_{\R^d \setminus \Omh} = 0 \b\}.
\end{displaymath}
Let $Q_h$ denote the P0 midpoint interpolation operator, so that
$\int_{\Omh} Q_h f$ is the midpoint rule approximation to $\int_{\Omh}
f$. 

Further, let $\beta \in C^{2,1}(\R^d)$ with $\beta = 0$ in $B_{\Ra}$
with $\Rdef \leq \Ra < \Rc$ and $\beta = 1$ in $\R^d \setminus
\Omh$, then we define the BQCE energy functional
\begin{align*}
  \Eb(u_h) := \sum_{a \in \L \cap \Omh} (1-\beta(a)) \Phi_a'(u) + \int_{\Omh}
  Q_h \b[ \beta W'(\D u_h) \b].
\end{align*}

\def\ub{u^{\rm b}_h}

The BQCE problem is to compute
\begin{equation}
  \label{eq:scheme:bqce}
  \ub \in \arg\min \b\{ \Eb(v_h) \bsep v_h \in \Ush \b\}.
\end{equation}

\subsection{Review of error estimates}
We present a formal review of error estimates for the BQCE scheme
established in \cite{2012-CMAME-optbqce,LiOrShVK:2014}. This
discussion will motivate our construction of the \BGFC scheme in the
next section.

Under a range of technical assumptions on $\Th$ and $\beta$ it is
shown in \cite{LiOrShVK:2014} that, if $\ua$ is a strongly stable
(positivity of the hessian) solution to \eqref{eq:atm}, and $(\beta,
\Th)$ are ``sufficiently well-adapted to $\ua$'', then there exists a
solution $\ub$ to \eqref{eq:scheme:bqce} such that
\begin{equation}
  \label{eq:scheme:bqce_errest}
  \| \D \ub - \D \overline{\ua} \|_{L^2} 
  \leq C_1 \| \D^2 \beta \|_{L^2} + C_2 \B( \|
  \beta h \D^2 \widetilde{\ua} \|_{L^2(\Omh)} +  \| \D \widetilde{\ua}
  \|_{L^2(\R^d\setminus B_{\Rc/2})} \B) + \dots,
\end{equation}
where ``$\dots$'' denotes formally higher order terms,
$\overline{\ua}$ denotes a P1 interpolant on the atomistic grid $\L$
and $\widetilde{\ua}$ a $C^{2,1}$-conforming interpolant on the
atomistic grid $\L$ (intuitively, $\D^j \widetilde{\ua}$, $j \geq 2$,
measure the local regularity of $\ua$; see \S~\ref{sec:analysis} and
\cite{LiOrShVK:2014} for more details). The constants depend on
(derivatives of) $\Phi_a, a \in \L,$ in a way that we will discuss in
more detail below.

The term $\| \beta h \D^2 \widetilde{\ua} \|_{L^2(\Omh)}$ measures the
finite element approximation error while the term $\| \D
\widetilde{\ua} \|_{L^2(\R^d\setminus B_{\Rc/2})}$ measures the error
committed by truncating to a finite computational domain. Exploiting
the generic decay rates \cite{EhrOrtSha:defectsV1} 
% \lz{LZ: reference [7] and [8] should be merged} 
% CO: done
\begin{equation}
  \label{eq:decay}
  |\D^j \widetilde{u^\a}(x)| \lesssim |x|^{1-d-j} \quad \text{as $j
    \to \infty$},
\end{equation}
these terms can be balanced by ensuring that $\Rc \approx (\Ra)^2$ and
the mesh is coarsened according to $h(x) \approx (|x|/\Ra)^{3/2}$ (see
\cite{2012-MATHCOMP-qce.pair}), which yields
\begin{displaymath}
  \| \beta h \D^2 \widetilde{\ua} \|_{L^2(\Omh)} +  \| \D \widetilde{\ua}
  \|_{L^2(\R^d\setminus B_{\Rc/2})} \lesssim (\dof)^{-1/2-1/d}.
\end{displaymath}

By contrast, the term $\|\D^2\beta\|_{L^2}$ is due to the (smeared)
ghost forces, and even an optimal choice of $\beta$ (balanced against
the atomistic region radius $\Ra$) yields only
\begin{equation}
  \label{eq:scheme:optimal_beta_rate}
  \|\D^2 \beta \|_{L^2} \approx (\dof)^{1/2-2/d}.
\end{equation}
To see this, we note that, a quasi-optimal choice is $\beta(x) =
B(r)$, where $B$ is a radial spline with $B(r) = 0$ for $r \leq \Ra$
and $B(r) = 1$ for $r \geq \Rb \in (\Ra, \Rc)$ (see
\cite{VkLusk2011:blended,2012-CMAME-optbqce} for in-depth
discussions).
% \lz{LZ: in-depth discussions?}). 
A straightforward computation (assuming $\Rb \lesssim
\Ra$; the case $\Rb \gg \Ra$ is similar) then shows that
$\|\D^2\beta\|_{L^2}^2 \approx (\Rb)^{d-1} (\Rb-\Ra)^{-3} +
(\Rb-\Ra)^{d-4}$, which is optimised subject to fixing the total
number of degrees of freedom if $\Rb - \Ra \approx \Ra$. This yields
precisely \eqref{eq:scheme:optimal_beta_rate}.

In summary, for this simple model problem, the BQCE scheme's rate of
convergence,
\begin{displaymath}
  \| \D \ub - \D \overline{\ua} \|_{L^2}  \lesssim (\dof)^{1/2-2/d},
\end{displaymath}
is the same in 2D and worse in 3D than a straightforward truncation
scheme, in which the atomistic model is minimised over a finite
computational domain (see~\cite{EhrOrtSha:defectsV1}). Analogous
results hold also for dislocations.

\subsection{The \BGFC scheme}
\label{sec:scheme:rbqce}
\def\Ebr{\E^{\rm bg}}
\def\ubr{u^{\rm bg}_h}
The motivation for the \BGFC scheme is to optimise the coefficient
$C_1$ in \eqref{eq:scheme:bqce_errest}. An investigation of the
analysis in \S~6.1 and \S~6.2 in \cite{LiOrShVK:2014} reveals that a
simple upper bound is
\begin{displaymath}
  C_1 \lesssim \sup_{a \in \L \cap {\rm supp}(\beta)} \sup_{b \in \L
    \setminus \{a\}}\B| \frac{\pp \Phi_a(y)}{\pp y(b)}|_{y = x+u} \B|.
\end{displaymath}
(As a matter of fact, this form of $C_1$ requires a minor modification
of the remaining error estimates \cite{LiOrShVK:2014}; however, we use
it only for motivation.)

The idea is to ``renormalise'' the interatomic potential so that $\del
\Phi_a(x) = 0$ for $|a|$ sufficiently large, which would then ensure
that $|\pp_{y(b)} \Phi_a(x+u)| \lesssim |u(b)-u(a)|$ and hence would
yield additional decay of the constant $C_1$ as the atomistic region
increases.

Recalling the discussion in \S~\ref{sec:scheme:hom} we note that
$\del\Phi_a''(x) = 0$ for all $a \in \L$. Thus, if we apply the
blending procedure to the renormalised atomistic energy
\eqref{eq:scheme:Ea_renormalised} then we would obtain a new constant
$C_1''$ with
\begin{align*}
  C_1'' &\lesssim \sup_{a \in \L \cap {\rm supp}(\nabla\beta)} \sup_{\substack{b \in \L
    \setminus \{a\} \\ r_{ab} \leq \rcut}} \b| \pp_{y(b)}
\Phi_a''(x+u) \b| \\
  &= \sup_{a \in \L \cap {\rm supp}(\nabla\beta)} \sup_{\substack{b \in \L
    \setminus \{a\} \\ r_{ab} \leq \rcut}} \b| \pp_{y(b)} \Phi_a(x+u)
- \pp_{y(b)} \Phi_a(x) \b| \\
&\lesssim C_2 \b\| \D \overline{\ua} \b\|_{L^\infty(\R^d \setminus
  B_{\Ra-2\rcut})} \lesssim (\Ra)^{-d}.
\end{align*}
where the second to last inequality is true for sufficiently large
$\Ra$ and the last inequality follows from the decay estimate given in
\cite[Thm. 3.1]{EhrOrtSha:defectsV1}. We therefore obtain
\begin{equation}
\label{eq:scheme:gfbound}
  C_1'' \|\D^2 \beta \|_{L^2} \lesssim (\dof)^{-1/2-2/d}, 
  % \lz{\text{LZ: should be $C_1''$ right?}} CO: yes
\end{equation}
which not only balances the best approximation error but is even
dominated by it.

To summarize, the \BGFC energy (difference) functional reads
\begin{equation}
  \label{eq:scheme:rbace}
  \Ebr(u_h) := \sum_{a \in \L \cap \Omh} (1-\beta(a)) \Phi_a''(u_h) +
  \int_{\Omh} Q_h\b[ \beta W''(\D u_h) \b] + \< \Lren, u_h \>,
\end{equation}
where $\Phi_a''$ is defined in \eqref{eq:scheme:defn_Phi''}, $\Lren$
in \eqref{eq:scheme:defn_Lren} and $W''(\mF) := W(\mI+\mF) - W(\mI) -
\pp W(\mI) : \mF$. The associated variational problem is
\begin{equation}
  \label{eq:scheme:rbqce}
  \ubr \in \arg\min \b\{ \Ebr(v_h) \bsep v_h \in \Ush \b\}.
\end{equation}

We can further optimise the \BGFC scheme as follows.
% \begin{remark}
%   \label{rem:scheme:rbqce_p1_balance}
If $\Rb / \Ra \sim c$ as $\Ra \to \infty$, then the coupling error of
the \BGFC scheme scales like $(\Ra)^{-d/2-2}$, and is therefore
dominated by the best approximation error, which scales like
$(\Ra)^{-d/2-1}$. To reduce computational cost (by a constant factor),
we can balance these two terms. Making the ansatz $\Rb - \Ra \sim
(\Ra)^t$, for $t \in (0, 1)$, and noting that we can always construct
$\beta$ such that $|\D^2\beta| \lesssim (\Ra)^{-t}$, we obtain that
\begin{displaymath}
   C_1'' \|\D^2 \beta \|_{L^2} \lesssim (\Ra)^{-d/2-3/2-1/2}.
\end{displaymath}
This is balanced with the best approximation rate, $(\Ra)^{-d/2-1}$,
if $t = 1/3$.

  % \begin{align*}
  %   \b[ C_2'' \|\D^2 \beta \|_{L^2}\b]^2 &\lesssim (\Ra)^{-2d} \b[
  %   (\Ra)^{d-1} (\Ra)^{-3t} + (\Ra)^{(d-4)t} \b] \\
  %   &\approx (\Ra)^{-d-1-3t} + (\Ra)^{-2d+t (d-4)}.\lz{\text{LZ: still, should be $C_1''$ right?}}
  % \end{align*}
% Thus, choosing $t \geq 1/3$ yields
%   \begin{displaymath}
%     C_1'' \|\D^2 \beta \|_{L^2} \lesssim (\Ra)^{-d-1} \approx (\dof)^{-1-1/d},  
%   \end{displaymath}
%   thus balancing the best approximation error \lz{LZ: this is rather dominated by the best approximation error if it is shown in the formula after (2.9)?}.
% % \end{remark}

We therefore conclude that, if $\Rb-\Ra \approx (\Ra)^t$ for some $t
\geq 1/3$, then we expect the \BGFC scheme to obey the error estimate
\begin{displaymath}
  \| \D \ubr - \D \overline{u^\a} \|_{L^2} \lesssim (\Ra)^{-d/2-1}
  \approx (\dof)^{-1/2-1/d}. 
  % \lz{\text{$(\dof)^{-1/2-1/d}$ and $(\Ra)^{-d/2-1}$?}}
\end{displaymath}
We shall make this rigorous for a slightly simplified formulation in
\S~\ref{sec:analysis}, where we will also prove an energy error
estimate.

\subsection{Connection to ghost-force correction and generalisation}
\label{sec:connection}
Consider, for simplicity, the case when $\Phi^{\rm hom}_a \equiv
\Phi_a$, i.e., the crystal is homogeneous. In this case, $L^{\rm ren}
\equiv 0$ as well. Moreover, we can rewrite the BGFC scheme as
follows:
\begin{align}
  \notag
  \Ebr(u_h) &= \Eb(u_h) - \sum_{a \in \L} (1-\beta(a)) \< \del
  \Phi_a(0), u_h \> - \int_{\R^d} Q_h\b[ \beta \pp W(0) : \D u_h \b]
  \dx \\
  \notag
  &= \Eb(u_h) - \< \del \Eb(0), u_h \> \\
  \label{eq:bgfc_gfform}
  &= \Eb(u_h) - \< \del \Eb(0) - \mathcal{F}^{\rm bqcf}(0), u_h \>,
\end{align}
where $\mathcal{F}^{\rm bqcf}$ is the BQCF operator defined in
\cite{2013-CMAME-bqcfcomp}, and $\mathcal{F}^{\rm bqcf}(0) = 0$ (this
non-conservative a/c coupling has no ghost forces). Thus, we see that
the renormalisation step $\Phi'_a \leadsto \Phi_a''$
(cf. \eqref{eq:scheme:Ea} and \eqref{eq:scheme:defn_Phi''}) 
% \lz{(LZ: it probably would be better
% to remind the reader the definition of $\Phi'_a$  in \eqref{eq:scheme:Ea} as it does not
% appear much)}
 is equivalent to
the dead load ghost-force correction scheme of Shenoy at al
\cite{Shenoy:1999a}, applied for a blended coupling formulation and in
the reference configuration.

This immediately suggests the following generalisation of the BGFC
scheme:
\begin{equation}
  \label{eq:bgfc_alternative}
  \Ebr(u_h) := \Eb(u_h) - \b\< \del\Eb(\hat{u}_h) - \mathcal{F}^{\rm
    bqcf}(\hat{u}_h), u_h - \hat{u}_h \b\>,
\end{equation}
where $\hat{u}_h$ is a suitable reference configuration, or
``predictor'', that can be cheaply obtained.

We will explore this alternative point of view in future work, in
particular with an eye to applications involving cracks and edge
dislocations.

% 2. If we choose $\Rb - \Ra \approx 1$ (e.g., only one atomic
%   spacing), then our method becomes a conservative variant of the
%  . In essence, instead of iterating the ghost
%   force correction to self-consistency, we only apply it in a
%   reference configuration. In this case, we would obtain that $C_1''
%   \|\D^2 \beta \|_{L^2} \lesssim (\Ra)^{-d/2-1/2} \approx
%   (\dof)^{-1/2-1/(2d)}$. By contrast, one expects the optimal rate
%   $(\dof)^{-1/2-1/d}$ for the force-based quasicontinuum method, which
%   is the self-consistent limit of the ghost force correction. There
%   is, at present, no proof of this fact, but for the related BQCF
%   scheme this rate is made rigorous in \cite{LiOrShVK:2014}.}

\section{Numerical Tests}
\label{sec:numerical_tests}
\subsection{Model problems}
Our prototype implementation of \BGFC is for the 2D triangular lattice
$\mA \Z^2$ defined by
\begin{displaymath}
  \mA = \begin{pmatrix} 1 & \cos(\pi/3) \\ 0 & \sin(\pi/3) \end{pmatrix}.
\end{displaymath}
To generate a defect, we remove $k$ atoms
\begin{displaymath}
  \cases{ \L^{\defc}_{k} := \b\{ - (k/2+1) e_1, \dots , k/2 e_1\b\}, &
    \text{ if $k$ is even, } \\[1mm]
    \L^{\defc}_{k} :=  \b\{ - (k-1)/2 e_1, \dots, (k-1)/2 e_1 \b\}, & \text{ if $k$ is odd,}
  }
\end{displaymath}
to obtain $\L := \mA\Z^2 \setminus \L^{\defc}_k$. For small $k$, the
defect acts like a point defect, while for large $k$ it acts like a
small crack embedded in the crystal. In our experiments we shall
consider $k = 2$ (di-vacancy) and $k = 11$ (microcrack), following
\cite{2012-CMAME-optbqce,2013-CMAME-bqcfcomp, OrZh:2013}.

The site energy is given by an EAM (toy-)model \eqref{eq:eam} \cite{Daw:1984a},
for which $\Phi_\ell$ is of the form 
\begin{align}
  \label{eq:eam}
  \Phi_\ell(y) &= \sum_{\rho \in \Rg(\ell)} \phi\b(|D_\rho y(\ell)|\b) + F\B(
  {\textstyle \sum_{\rho \in \Rg(\ell)}} \psi\b( |D_\rho y(\ell)|\b)
  \B), \\
  \notag
  \text{with} \qquad &\phi(r) = [e^{-2a(r-1)}-2e^{-a(r-1)}], \quad \psi(r) = e^{-br}, \\
  \notag
  &F(\tilde{\rho})=c\b[(\tilde{\rho}-\tilde{\rho}_0)^2+(\tilde{\rho}-\tilde{\rho}_0)^4\b],
\end{align}
with parameters $a = 4.4, b = 3, c=5, \tilde{\rho}_0 = 6e^{-b}$. The
interaction range is $\Rg(\ell) = \L \cap B_2(\ell)$, i.e., next
nearest neighbors in hopping distance.

To construct the BQCE and \BGFC schemes, we choose an elongated
hexagonal domain $\Omega^\a$ containing $K$ layers of atoms
surrounding the vacancy sites and the full computational domain
$\Omega_h$ to be an elongated hexagon containing $\Rc$ layers of atoms
surrounding the vacancy sites. The domain parameters are chosen so
that $\Rc = \lceil \smfrac12 (\Ra)^2 \rceil$. The finite element mesh
is graded so that the mesh size function $h(x) = {\rm diam}(T)$ for $T
\in \TT$ satisfies $h(x) \approx (|x|/\Ra)^{3/2}$. These choices balance
the coupling error at the interface, the finite element interpolation
error and the far-field truncation error
\cite[Sec.5.2]{EhrOrtSha:defectsV1}. Recall moreover that $\dof :=
\#\Nh$.

The blending function is obtained in a preprocessing step by
approximately minimising $\|\D^2\beta\|_{L^2}$, as described in detail
in \cite{2012-CMAME-optbqce}.

We implement the equivalent ghost force removal formulation
\eqref{eq:bgfc_gfform} instead of the ``renormalisation formulation''
\eqref{eq:scheme:rbace}.

\subsubsection{Di-vacancy}
\label{sec:di-vacancy}
In the di-vacancy test two neighboring sites are removed, i.e., $k =
2$. We apply $3\%$ isotropic stretch and $3\%$ shear loading, by
setting
\begin{displaymath}
\mB := \begin{pmatrix} 1+s & \gamma_{\rm II} \\ 0 & 1+s \end{pmatrix}\cdot \mF_0.
\end{displaymath}
where $\mF_0 \propto I$ minimizes $W$, $s=\gamma_{\rm II}=0.03$.

% One then obtains
%   \cite[Prop. 5.5]{EhrOrtSha:defectsV1} under additional conditions on the
%   stability of the method and the magnitude of the reconstruction
%   parameters (we can verify both only {\it a posteriori}) that
% \begin{equation}
%   \label{eq:errest}
%   \| \D y - \D y_h \|_{L^2} \leq C {\rm DOF}^{-1},
% \end{equation}
% where $y$ is identified with its P1 interpolant on the canonical
% triangulation of $\L$ and ${\rm DOF}$ denotes the total number of
% degrees of freedom (i.e. the number of atomistic sites $\L^{\a,\i}$
% plus the number of finite element nodes).

% The error in the $H^1$-seminorm is displayed in Figure
% \ref{fig:rate_divacancy_w12}, the errors in the
% $W^{1,\infty}$-seminorm are displayed in Figure
% \ref{fig:rate_divacancy_w1inf}, the errors in the energy are displayed
% in Figure \ref{fig:rate_divacancy_E}. For the $H^1$ and
% $W^{1,\infty}$-seminorms we now observe precisely the rate of
% convergence $\dof^{-1/2}$ and $\dof^{-1}$.

\subsubsection{Microcrack}
\label{sec:micro-crack}

In the microcrack experiment, we remove a longer segment of atoms,
$\L^{\defc}_{11}=\{-5e_1,\dots,5e_1\}$ from the computational
domain. The body is then loaded in mixed mode ${\rm I}$ \& ${\rm II}$,
by setting,
\begin{displaymath}
\mB := \begin{pmatrix} 1 & \gamma_{\rm II} \\ 0 & 1+\gamma_{\rm I} \end{pmatrix}\cdot \mF_0.
\end{displaymath}
where $\mF_0\propto I$ minimizes $W$, and $\gamma_{\rm I}=\gamma_{\rm
  II}=0.03$ ($3\%$ shear and $3\%$ tensile stretch).

\subsection{Methods}
We shall test \BGFC method with blending widths $K := \Rb-\Ra= \lceil
(\Ra)^{1/3} \rceil$ and with $K = \Rb-\Ra = \Ra$ (here, the blending
width denotes the number of hexagonal atomic layers in the blending
region).  The \BGFC scheme is compared against the 3 competitors
previously considered in \cite{2012-CMAME-optbqce,
  2013-CMAME-bqcfcomp, OrZh:2013}:
\begin{itemize}
\item {\bf B-QCE}: blended quasicontinuum method,
  implementation based on \cite{2012-CMAME-optbqce}, with most details
  described in \S~\ref{sec:scheme:bqce}. 
  % ; B-QCE+ is a variant
  % with highly optimised approximation parameters described in
  % \cite[Sec. 4.3]{2013-CMAME-bqcfcomp}.
\item {\bf GRAC}: sharp-interface consistent energy-based a/c coupling
  \cite{OrZh:2013}.
\item {\bf B-QCF}: blended force-based a/c coupling, as described in
  \cite{2013-CMAME-bqcfcomp}. Energies of B-QCF are computed using
  B-QCE (i.e., the B-QCE energy is evaluated at the B-QCF solution).
\end{itemize}

\subsection{Results}
We
present two experiments, a di-vacancy ($k = 2$) and a microcrack
($k = 11$). For each test, we choose an increasing sequence of
atomistic region sizes $\Ra$, followed by the quasi-optimal choices of
$\Rb, \Omega_h, \beta$ as described above.

For both experiments we plot the absolute errors against the number of
degrees of freedom (DOF), which is proportional to computational cost,
in the $H^1$-seminorm, the $W^{1,\infty}$-seminorm and in the
(relative) energy.

\begin{figure}
  \begin{center}
    \includegraphics[width=12.5cm]{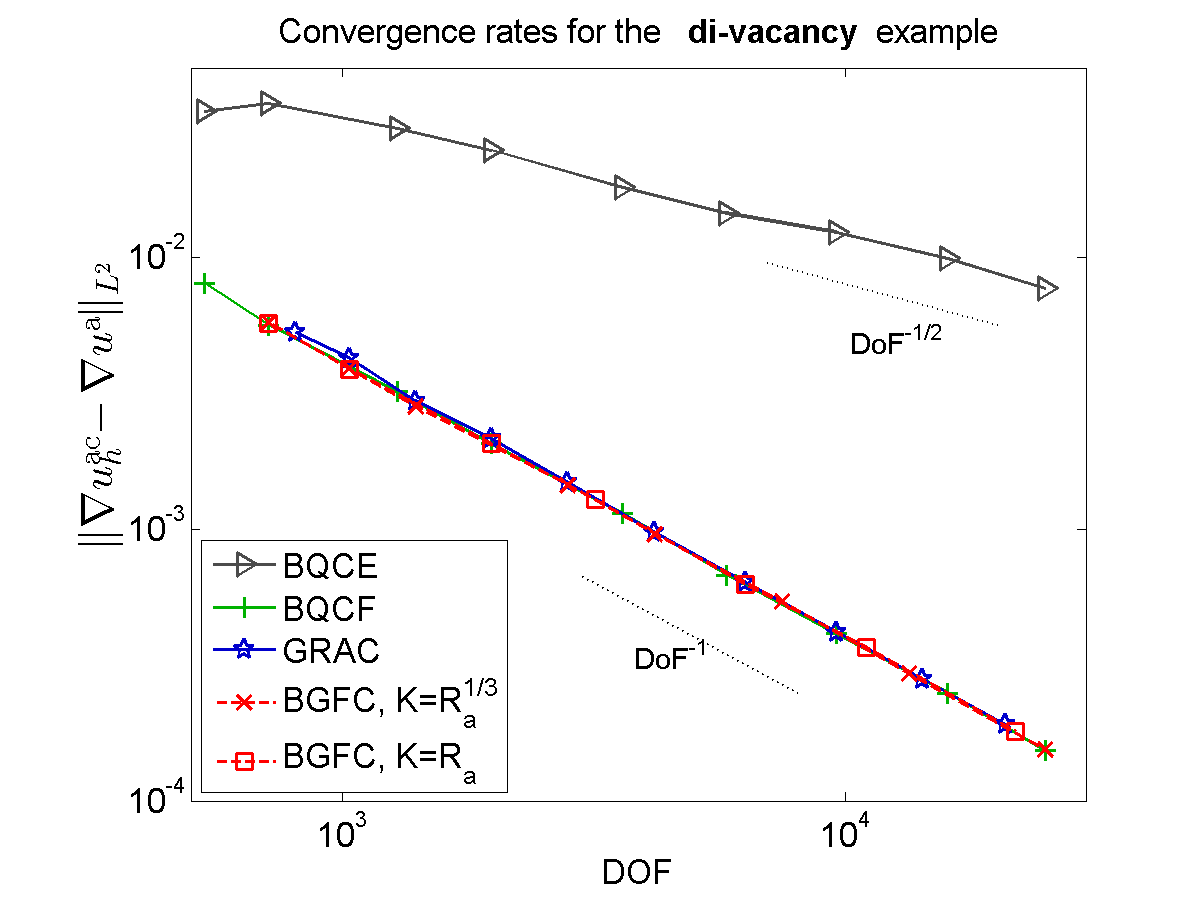}
    \caption{\label{fig:rate_divacancy_w12} Convergence rates in the
      energy-norm (the $H^1$-seminorm) for the di-vacancy benchmark
      problem described in Section \S~\ref{sec:di-vacancy} . }
  \end{center}
\end{figure}

\begin{figure}
  \begin{center}
    \includegraphics[width=12.5cm]{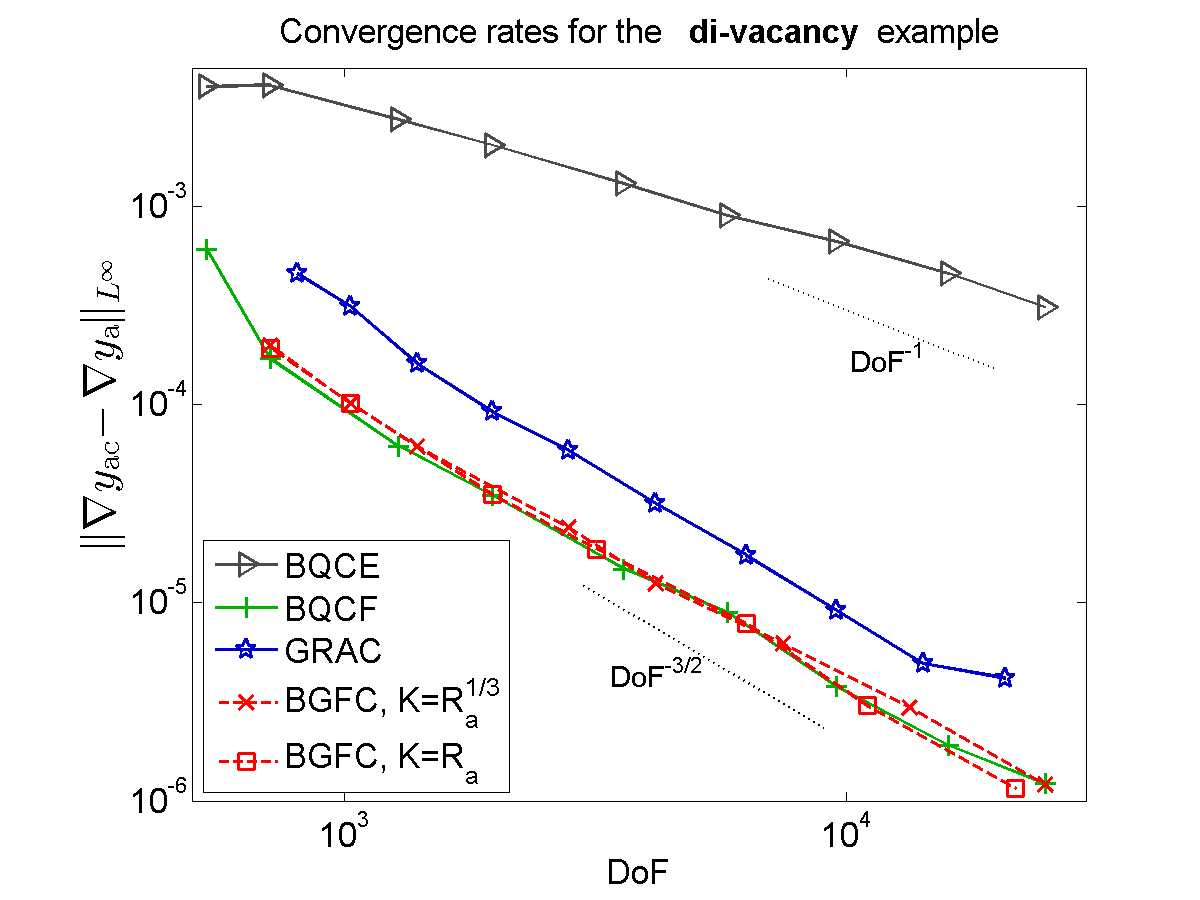}
    \caption{\label{fig:rate_divacancy_w1inf} Convergence rates in the $W^{1,\infty}$-seminorm for the di-vacancy benchmark problem described in Section \S~\ref{sec:di-vacancy} .}
  \end{center}
\end{figure}

\begin{figure}
  \begin{center}
    \includegraphics[width=12cm]{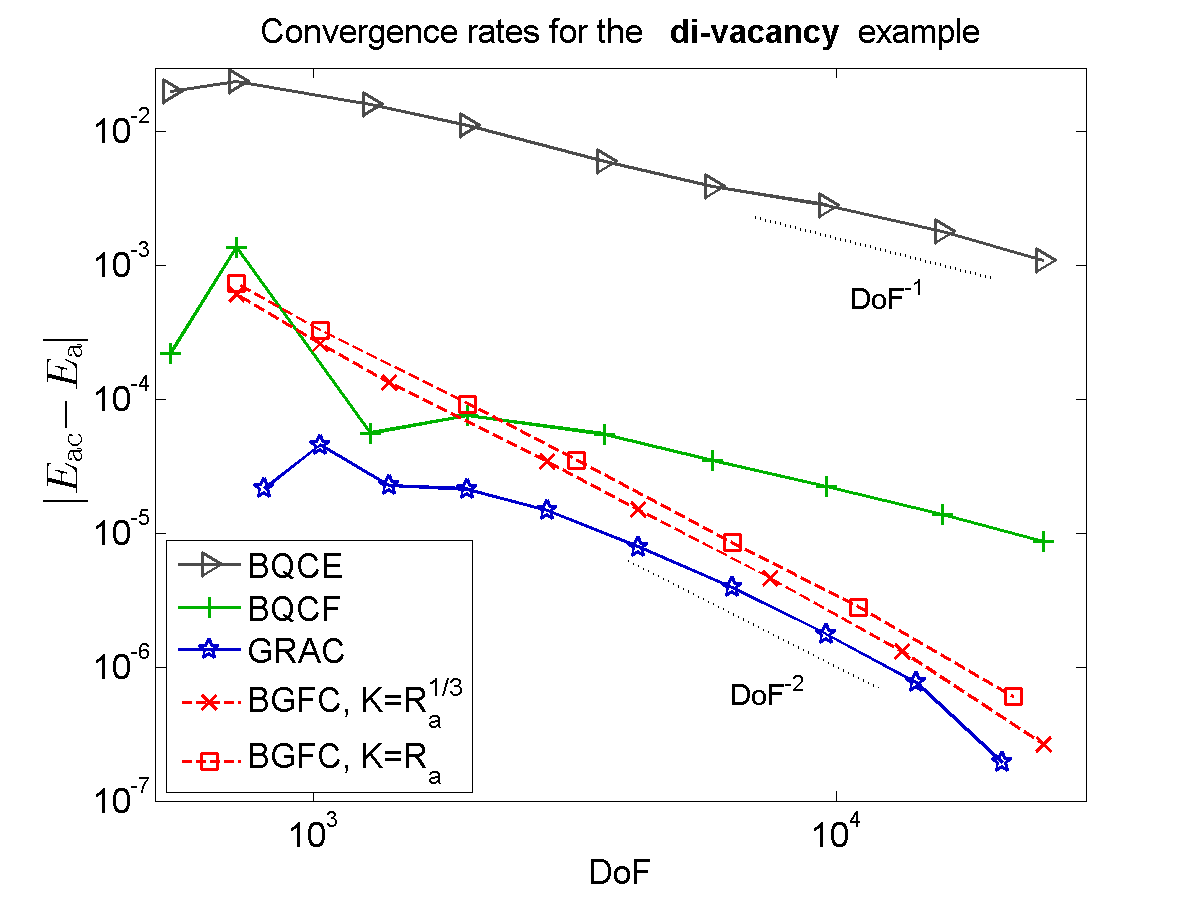}
    \caption{\label{fig:rate_divacancy_E} Convergence rates in the
      {relative energy} for the di-vacancy benchmark problem described
      in Section \S~\ref{sec:di-vacancy} .}
  \end{center}
\end{figure}

\begin{figure}
  \begin{center}
    \includegraphics[width=12cm]{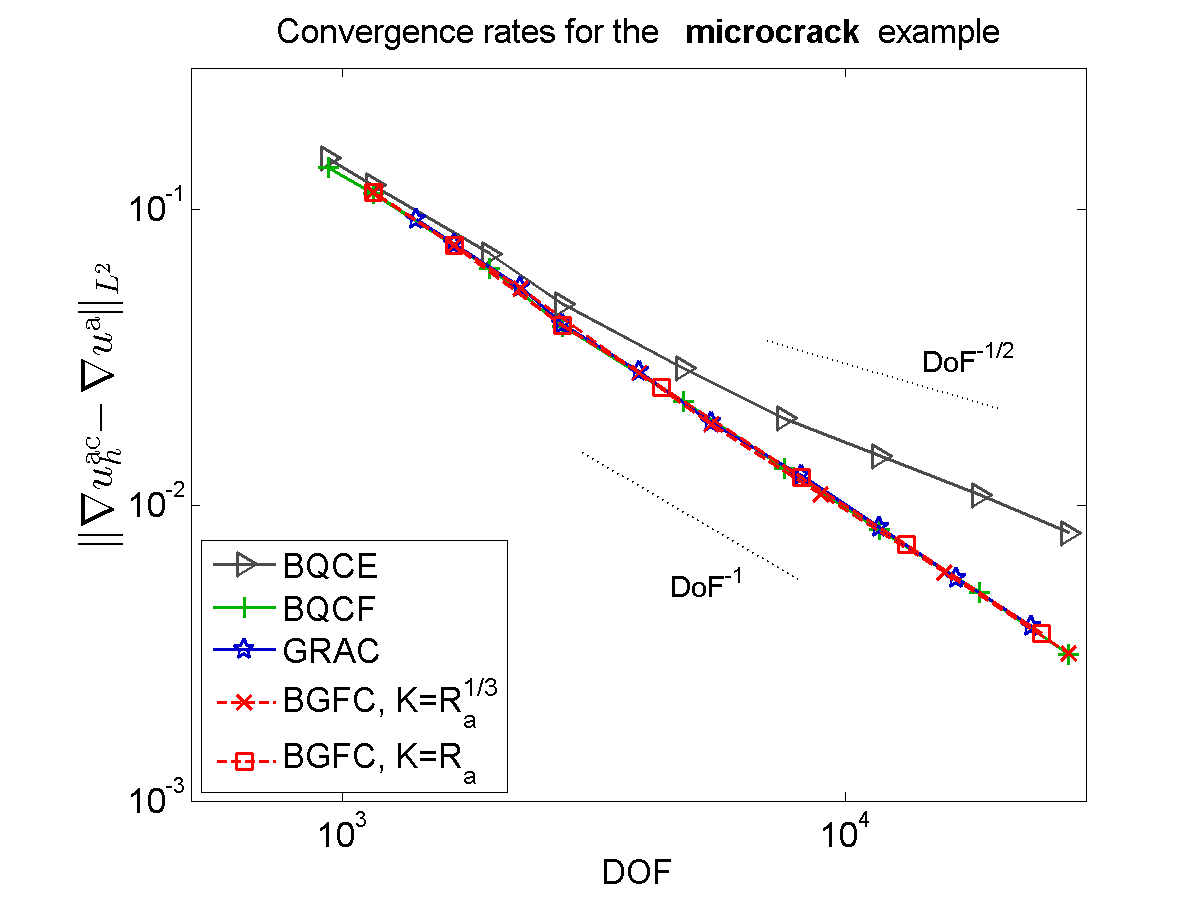}
    \caption{\label{fig:rate_mcrack_w12} Convergence rates in the energy-norm (the $H^1$-seminorm) for the microcrack benchmark problem with $\L^{\defc}_{11}$ described in Section \S~\ref{sec:micro-crack} .}
  \end{center}
\end{figure}

\begin{figure}
  \begin{center}
    \includegraphics[width=12cm]{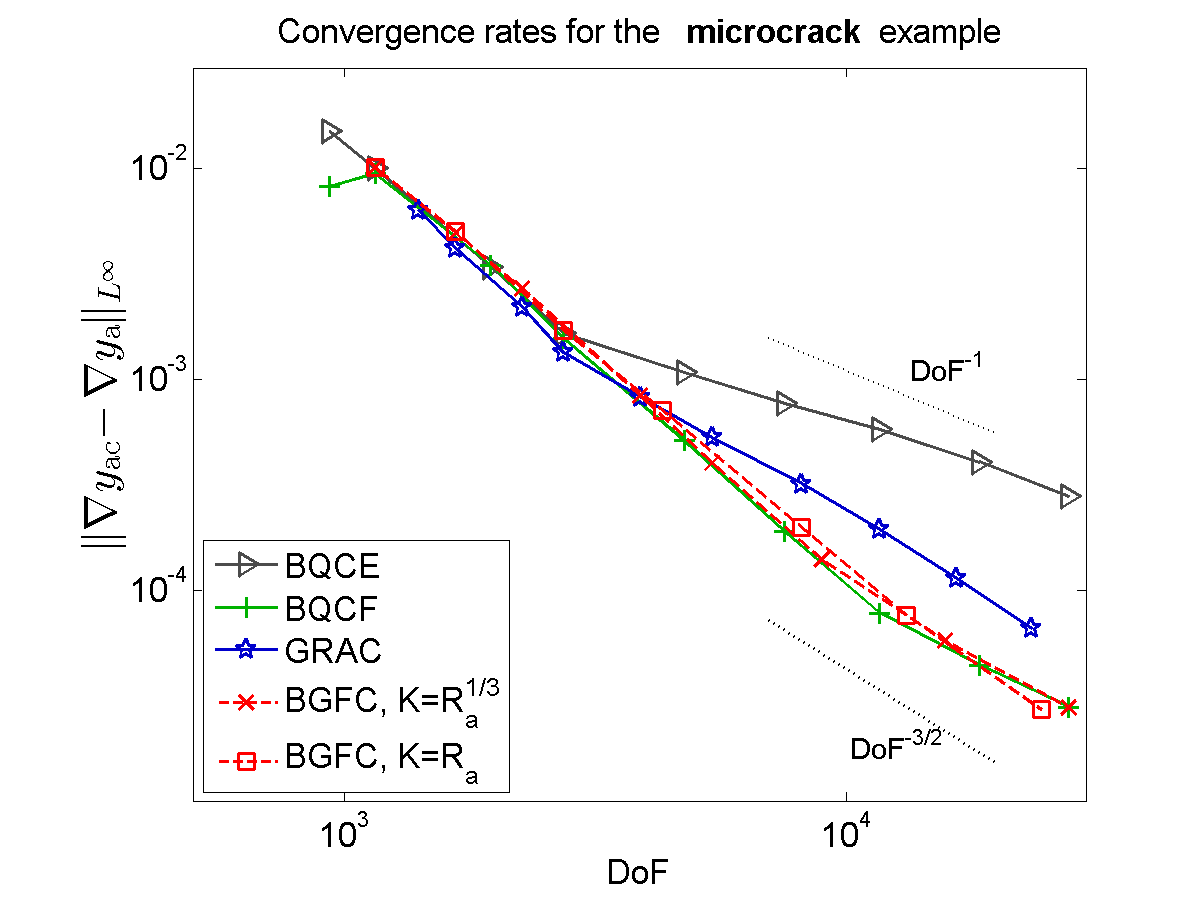}
    \caption{\label{fig:rate_mcrack_w1inf} Convergence rates in
      the {$W^{1,\infty}$-seminorm} for the microcrack benchmark
      problem with $\L^{\defc}_{11}$ described in Section
      \S~\ref{sec:micro-crack} .}
  \end{center}
\end{figure}

\begin{figure}
  \begin{center}
    \includegraphics[width=12cm]{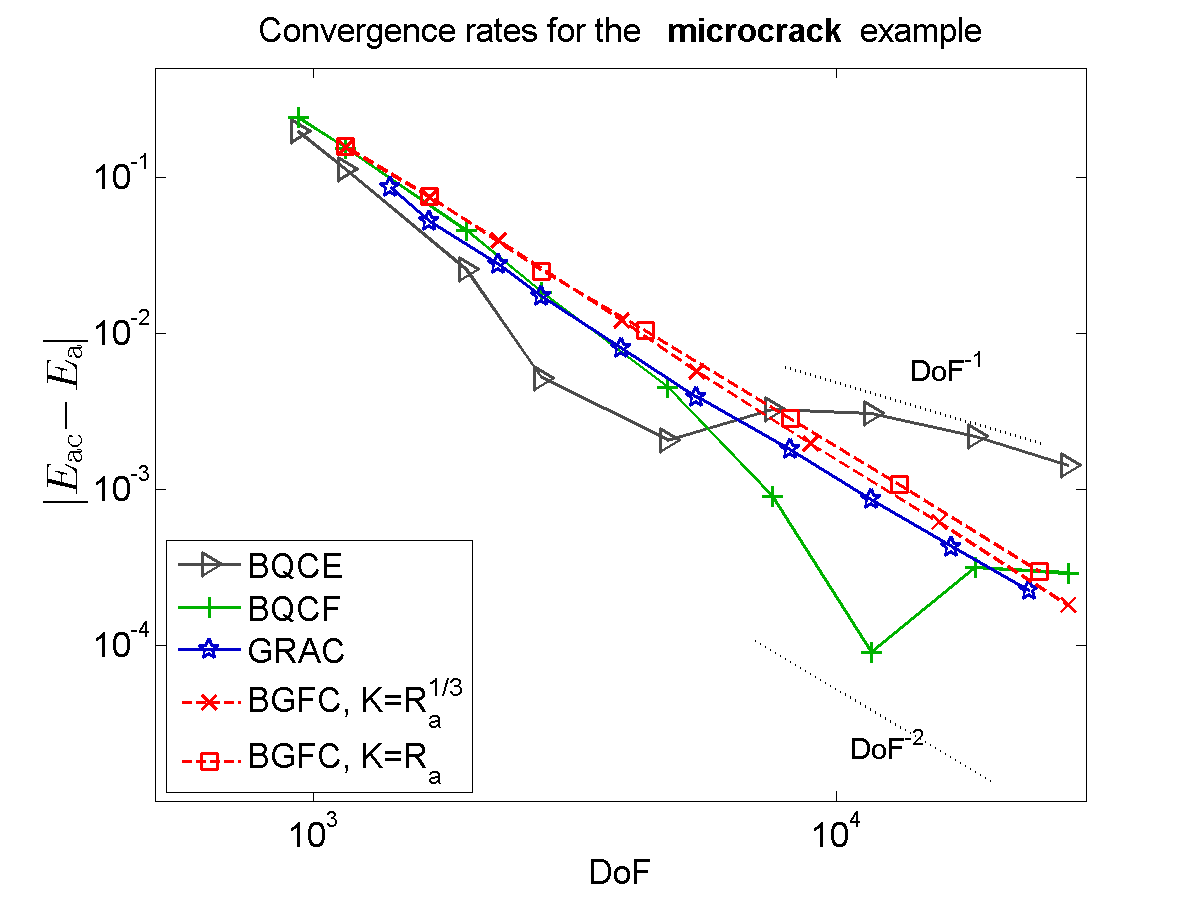}
    \caption{\label{fig:rate_mcrack_E} Convergence rates in the
      {relative energy} for the microcrack benchmark problem with
      $\L^{\defc}_{11}$ described in Section \S~\ref{sec:micro-crack}
      .}
  \end{center}
\end{figure}

The results are shown in Figures \ref{fig:rate_divacancy_w12},
\ref{fig:rate_divacancy_w1inf} and \ref{fig:rate_divacancy_E} for the
di-vacancy problem and in Figures \ref{fig:rate_mcrack_w12},
\ref{fig:rate_mcrack_w1inf} and \ref{fig:rate_mcrack_E}
for the microcrack problem.

In the first experiment, we are able to clearly observe the predicted
asymptotic behaviour of the a/c coupling schemes, while in the second
experiment we observe a significant pre-asymptotic regime where the
analytic predictions become relevant only at fairly high resolutions.

In all error graphs we clearly observe the optimal convergence rate of
\BGFC, together with other consistent methods GRAC and BQCF, while
BQCE has a sub-optimal rate.

%Indeed, the errors are competitive with the quasi-optimal force-based schemes (QCF, B-QCF):
%for $H^1$ errors they are essentially comparable, for
%$W^{1,\infty}$ errors the force-based schemes are only better by a
%moderate constant factor, while for the energy errors the GRAC methods
%are optimal. (Note that, for QCF we evaluate the QCE energy and for
%B-QCF we evaluate the B-QCE energy.)

\section{Extensions}
\label{sec:scheme:ext}
It is possible to extend the formulation of the \BGFC scheme to a much
wider range of problems, including e.g. multiple defect regions,
problems with surfaces (e.g., nano-indentation, crack propagation),
complex crystals, or higher order finite elements. We now present a
range of such generalisations, arguing only formally to motivate a
more complete and rigorous treatment in future work.

\subsection{Higher-order finite elements}
\label{sec:ext:p2}
We have seen in \S~\ref{sec:scheme:rbqce} that in the BGFC scheme
applied to point defects, the approximation error dominates the
blending (coupling) error. This particular bottleneck is relatively
straightforward to remove by increasing the order of the finite
element scheme and the size of the continuum region. The following
discussion is motivated by \cite{WuP2}.

We construct the computational domain and finite element mesh in the
same way as in \S~\ref{sec:scheme:rbqce} and
\S~\ref{sec:numerical_tests}. We decompose $\Th = \Th^{({\rm P1})}
\cup \Th^{({\rm P2})}$, where
\begin{displaymath}
  \Th^{({\rm P1})} = \b\{ T \in \Th \bsep \beta|_T < 1 \b\},
\end{displaymath}
and replace $\Us_h$ with the
approximation space
\begin{align*}
  \Us_h^{(2)} := \b\{ u_h \in C(\R^d; \R^d) \bsep~& u_h |_T \text{ is
    affine for $T \in \Th^{({\rm P1})}$, and } \\
    &u_h|_T \text{ is quadratic for $T \in \Th^{({\rm P2})}$} \b\}.
\end{align*}
That is, we retain the P1 discretisation in the fully refined
atomistic and blending region where $\L$ and $\Nh$ coincide, but
employ P2 finite elements in the continuum region. Accordingly, the
qudrature operator $Q_h$ (previously midpoint interpolation) must now
be adjusted to provide a second-order quadrature scheme so that
$\nabla u_h \otimes \nabla u_h$ for $u_h \in \Us_h^{(2)}$ can be
integrated exactly.

The resulting P2-BGFC method reads
\begin{displaymath}
  u_h \in \arg\min \b\{ \Ebr(u_h) \bsep u_h \in \Us_h^{(2)} \b\}.
\end{displaymath}

\subsubsection{Convergence rate}
The blending error and the Cauchy--Born modelling error contributions
to the P2-BGFC method remain the same as for the P1-BGFC method,
$C_1'' \|\D^2\beta\|_{L^2} \lesssim (\Ra)^{-d/2-2}$; see
\S~\ref{sec:scheme:rbqce}, Equation \eqref{eq:scheme:gfbound}. Only
the approximation error component must be reconsidered. This requires
some non-trivial modifications to the analysis we present in
\S~\ref{sec:analysis}, which are left to future work, however it is
reasonable to expect that the best approximation error contribution
can be bounded by
\begin{displaymath}
  \| h^2 \D^3 \tilde{u}^\a \|_{L^2(\Omega_h \setminus B_{\Ra})} + \|
  \D \tilde{u}^\a \|_{L^2(\R^d \setminus B_{\Rc/2})},
\end{displaymath}
where the first term is the standard P2 finite element best
approximation error and the second term is the far-field truncation
error.

Choosing $h(x) \approx \max(1, (|x|/\Rb)^{3/2})$ and an increased
continuum region $\Rc \approx (\Ra)^{1+4/d}$, a straightforward
computation, employing the generic decay rates \eqref{eq:decay} for
point defects, shows that the two terms are balanced and bounded by
% \lesssim (\Ra)^{-3/2} \bg( \int_{\Ra}^{\infty} r^{d-1+3-2d-4} \dr
% \bg)^{1/2} \lesssim (\Ra)^{-3/2-d/2-1/2} 
% 
% \lesssim ( \int_{Rc}^\infty r^{d-1-2d})^{1/2} \approx (\Rc)^{-d/2}
% 
% td/2 = d/2+2  <>   (t-1) = 2 / (d/2) = 4/d   <>   t = 1 + 4/d.
%
\begin{align*}
  \| h^2 \D^3 \tilde{u}^\a \|_{L^2(\Omega_h \setminus B_{\Ra})} + \|
  \D \tilde{u}^\a \|_{L^2(\R^d \setminus B_{\Rc/2})}
  &\lesssim \bg( \int_{\Ra}^\infty r^{d-1} \B( \smfrac{r}{\Ra}\B)^3
  r^{-2d-4} \dr \bg)^{1/2} + (\Rc)^{-d/2} \\
  &\lesssim
 (\Ra)^{-d/2-2} \approx (\dof)^{-1/2-2/d}.
\end{align*}
Note that it is possible to make this construction without violating
the necessary mesh regularity required to obtain the stated finite
element approximation error; see \cite{OlBoLuSh:2013} for further
discussion.
%\lz{LZ: this sounds a bit abrupt, maybe explain more?}

Thus we (formally) obtain
\begin{displaymath}
  \| \D \ubr - \D \bar{u}^\a  \|_{L^2} \lesssim (\Ra)^{-d/2-2} 
  \approx (\dof)^{-1/2-2/d}.
\end{displaymath}

It is particularly interesting to note that the Cauchy--Born modelling
error contribution is also bounded by
\begin{displaymath}
  \| \D^3 \tilde{u}^\a \|_{L^2(\Omega_h \setminus B_{\Ra})} + \| \D^2
  \tilde{u}^\a \|_{L^4(\Omega_h \setminus B_{\Ra})}^2 \lesssim
  (\Ra)^{-d/2-2} \approx (\dof)^{-1/2-2/d}
\end{displaymath}
and this bound is in general optimal. Thus, we see that for the P2-GFC
method all three error components (coupling error, approximation
error, Cauchy--Born modelling error) are balanced in the
energy-norm. In particular, this means that, for point defects, the
P2-BGFC scheme is {\em quasi-optimal} among all a/c coupling method
that employ the Cauchy--Born model in the continuum region. We plan,
in future work, to present a complete analysis and implementation of
this scheme.

\subsection{Screw dislocation}
\label{sec:ext:disl}
We briefly demonstrate how the \BGFC scheme may be formulated for
simulating a screw dislocation. The ideas are a straightforward
combination of those in \cite{EhrOrtSha:defectsV1} and
\S~\ref{sec:scheme} of the present paper, thus we only present minimal
details.

For the sake of simplicity, we restrict the discussion and
implementation to nearest-neighbour interaction and anti-plane shear
motion, following \cite[Sec. 6.2]{EhrOrtSha:defectsV1}. That is, we
define $\L = \Lhom = \mA \Z^2$ where
\begin{displaymath}
  \mA = \left(\begin{matrix}
      1 & \cos(\frac{\pi}{3}) \\ 0 & \sin(\frac{\pi}{3}) \end{matrix}
  \right) \quad \text{and} \quad
  \Rg = \{ \mQ_6^j e_1 \sep j = 0, \dots, 5 \}
\end{displaymath}
is the set of interacting lattice directions.  (Note that $\Lhom$ is
in fact the projection of a bcc crystal along the (111) direction.)
Admissible (anti-plane) displacements are maps $y : \L \to \R$.

The site potential is now a map $\Phi_a \in C^3(\R^6)$, i.e.,
$\Phi_a(y)$ is a function of $(y(b)-y(a))_{b \in a+\Rg}$. To admit
slip by a Burgers vector (we assume the Burgers vector is $(1,0,0)$),
we assume that $\Phi_a(y) = \Phi_a(z)$ whenever $y - z : \L \to \Z$.

A screw dislocation is enforced (e.g.) by applying the far-field
boundary condition
\begin{displaymath}
  y(a) \sim y^{\rm lin}(a) := \smfrac{1}{2\pi} \arg(a - \hat{a}), \quad
  \text{as } |a| \to \infty,
\end{displaymath}
where $y^{\rm lin}$ is the linearised elasticity solution and
$\hat{a}$ an arbitrary shift of the dislocation core. The model of
\S~\ref{sec:scheme} can be extended by defining $\Phi_a'(u) :=
\Phi_a(y^{\rm lin}+u) - \Phi_a(y^{\rm lin})$, and $\Ea(u) := \sum_{a
  \in \L} \Phi_a'(u)$. With this modification the exact model still
reads \eqref{eq:atm}; see \cite{EhrOrtSha:defectsV1} for the details.

To define the \BGFC scheme we renormalise $\Phi_a$ a second time,
\begin{displaymath}
  \Phi_a''(u) := \Phi_a(y^{\rm lin}+u) - \Phi_a(y^{\rm lin}) - \< \del
  \Phi_a(0), u \>,
\end{displaymath}
which gives rise to the \BGFC functional defined by
\eqref{eq:scheme:rbace}. Note that $L^{\rm ren} \equiv 0$ in this
case. The resulting \BGFC scheme is still given by
\eqref{eq:scheme:rbqce}.

\begin{remark}
  It is tempting to define
  \begin{displaymath}
    \Phi_a''(u) := \Phi_a(y^{\rm lin}+u) - \Phi_a(y^{\rm lin}) - \<
    \del\Phi_a(y^{\rm lin}), u \>,
  \end{displaymath}
  which has seemingly has advantages in terms of error reduction.
  However, (i) it has the downside of having to evaluate a
  non-trivial functional $\< L^{\rm ren}, u \>$, for which a new
  scheme must be developed, and (ii) the Cauchy--Born modelling error
  is already dominant in the dislocation case, which means that no
  further improvements can in fact be expected.

  However, taking the alternative view presented in \S
  \ref{sec:connection} we may (re-)define the BGFC scheme as in
  \eqref{eq:bgfc_alternative}, where we note that now the predictor
  $y^{\rm lin}$ is used for the dead load ghost force removal {\em
    without} creating long-ranging residual forces in the continuum
  region.  This may indeed lead to a (moderate) improvement, which we
  will analyze in future work, together with applications to edge
  dislocations, where the ``renormalisation formulation'' seems less
  straightforward.
\end{remark}

\subsubsection{Convergence rate}
Suppose that the setup of the computational geometry is as in
\S~\ref{sec:scheme}, with the only exception that we now need to take
a coarsening rate $h(x) \approx |x| / \Ra$, due to the quadrature error (see
\cite{LiOrShVK:2014}). This only marginally modifies the analysis.

It is shown in \cite[Thm. 3.1]{EhrOrtSha:defectsV1} under natural
technical conditions that, if a minimiser $\ua$ of the screw
dislocation problem exists, then
\begin{displaymath}
  |\D^j y_0(x)| \lesssim |x|^{2-j} \quad \text{and} \quad
  |\D^j \widetilde{\ua}(x)| \lesssim |x|^{-1-j} \log |x|.
\end{displaymath}

A QNL type ghost-force free scheme (e.g., geometric-reconstruction
based \cite{OrZh:2013}) is then expected to have an error of the order
of magnitude of (see
\cite[Sec. 5.2]{EhrOrtSha:defectsV1} for the details)
\begin{align*}
  \| \D \overline{\ua} - \D u_h^{\rm qnl} \|_{L^2} &\lesssim \| h \D^2 \widetilde{\ua}
  \|_{L^2(\Omh \setminus B_{\Ra})} + \| \D \widetilde{\ua} 
  \|_{L^2(\R^d \setminus B_{\Rc / 2})} \\
  & \qquad + \| \D^2 (y_0+\widetilde{\ua}) \|_{L^2(\Omega^{\rm i})} + \dots \\
  &\lesssim (\Ra)^{-2} (\log \Ra)^{1/2} + (\Ra)^{-3/2}
  \approx (\dof)^{-3/4}.
\end{align*}
where $u_h^{\rm qnl}$ denotes the solution of such a scheme, $\| \D^2
(y_0+\widetilde{\ua}) \|_{L^2(\Omega^{\rm i})}$ the interfacial
coupling error (cf. \cite{Or:2011a}), and we denoted again several
dominated terms by ``$\dots$''. We observe that, for QNL type methods,
the coupling error dominates the estimate. 

Following the analysis in \cite{LiOrShVK:2014},
\S~\ref{sec:scheme:rbqce} and \S~\ref{sec:analysis}, we can obtain
\begin{align*}
    \| \D \overline{\ua} - \D \ubr \|_{L^2} &\lesssim \| \D^2\beta
    \|_{L^\infty} \| \D (y_0+\widetilde{\ua}) \|_{L^2(\Omb)} + {\rm best~approx.~err.} + \dots \\
    &\lesssim (\Rb-\Ra)^{-2} \b(\log\b( \Rb/\Ra \b)\b)^{1/2} + (\Ra)^{-2} (\log\Ra)^{1/2}.
\end{align*}
We observe that with $\Rb-\Ra \approx (\Ra)^{\alpha}$ the two errors
are balanced for $\alpha = 1$, i.e., $\Rb-\Ra \approx \Ra$, and that
in this case we obtain the ``optimal rate'' (i.e., the best
approximation error rate)
\begin{displaymath}
  \| \D \overline{\ua} - \D \ubr \|_{L^2} \lesssim (\Ra)^{-2}
  (\log\Ra)^{1/2} \approx (\dof)^{-1} (\log\dof)^{1/2}.
\end{displaymath}

Thus, we conclude that the \BGFC scheme leads to a better rate of
convergence than the QNL type scheme. This is particularly encouraging
as the latter is often assumed optimal among energy-based a/c coupling
schemes.

\begin{remark}
  We note that the Cauchy--Born modelling error for the screw
  dislocation example is bounded, in terms of $y = y_0+u$, by
  \begin{displaymath}
    \| \D^3 \tilde{y} \|_{L^2(\Omega_h \setminus B_{\Ra})} + \| \D^2
    \tilde{y} \|_{L^4(\Omega_h\setminus B_{\Ra})}^2 \lesssim
    (\Ra)^{-2} \approx (\dof)^{-1}.
  \end{displaymath}
  Thus, up to logarithmic terms, the best approximation and blending
  errors are both balanced with the Cauchy--Born modelling error,
  which is a lower-bound for a/c couplings based on local continuum
  models. In this sense, the \BGFC method is {\em optimal} for screw
  dislocations as well.
\end{remark}

\subsubsection{Numerical experiment}
Replicating the setting from \cite[Sec. 6.2]{EhrOrtSha:defectsV1}, we
use a simplified EAM-type interatomic potential
(cf. \eqref{eq:eam_potential}), given by
\begin{align*}
  \Phi_a(y) := G\B( \sum_{b \in a+\Rg} \phi(y(b)-y(a)) \B), \quad
  \text{where~}& \quad G(s) = 1+\smfrac12 s^2, \\[-3mm]
  \quad \text{and~}& \quad
  \phi(r) = \sin^2\b( \pi r \b).
\end{align*}
Note that, in this case, the BQCE and BGFC methods are in fact
identical since $\del \Phi_a(0) = 0$. (This is an artefact of the
anti-plane setting.) 
% \lz{LZ: is this also true for longer range interactions?}
% CO: for any finite interaction range

We employ the same constructions of the computational domain as in the
point defect case described in \S~\ref{sec:numerical_tests}, but
without vacancy sites.

\begin{figure}
  \begin{center}
    \includegraphics[height=8cm]{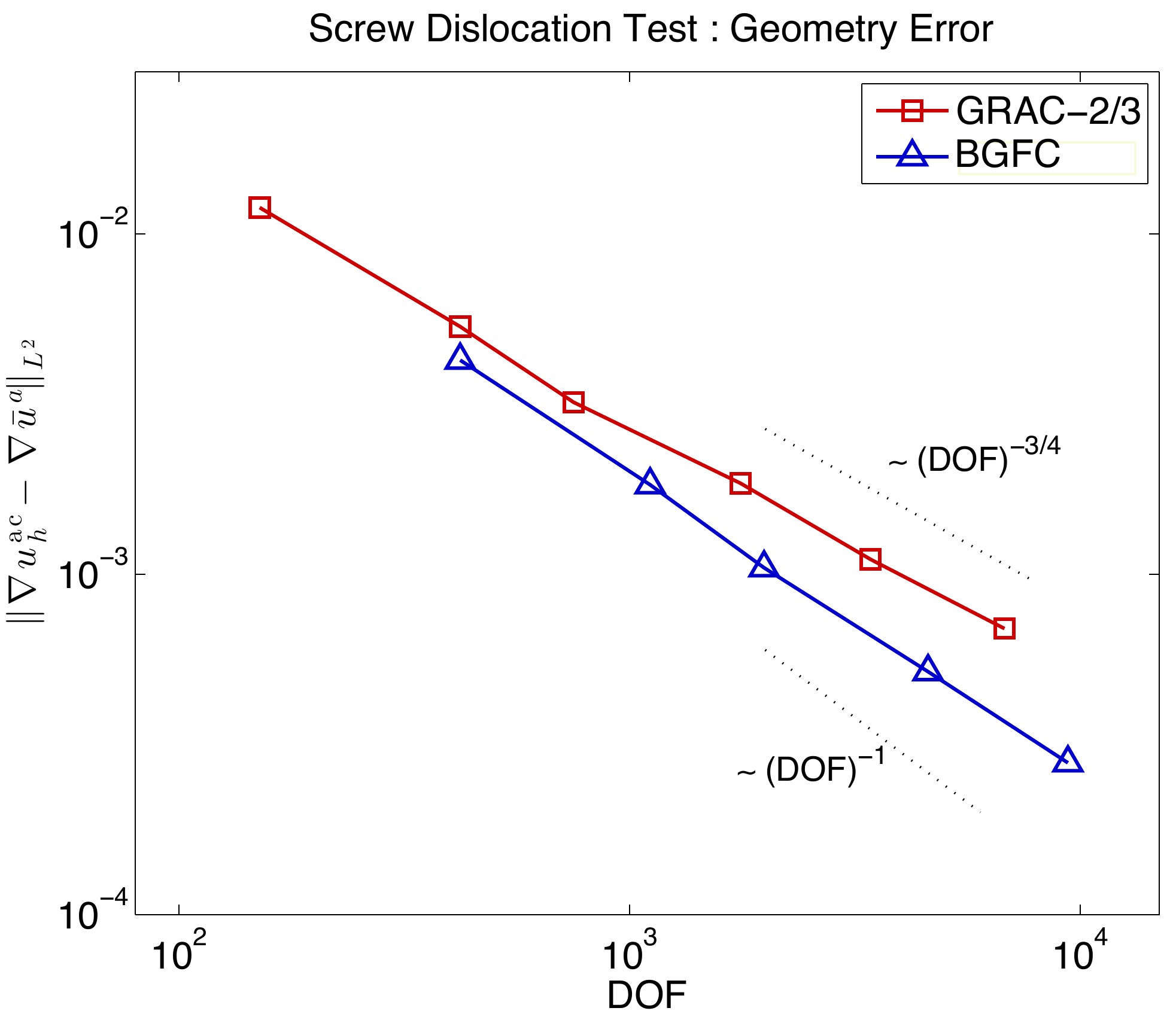}
    \caption{\label{fig:screwrb_err2} Convergence of the QNL-type
      GRAC-2/3 \cite{2012-OrtZha-SINUM} and of the \BGFC schemes for
      an anti-plane screw dislocation problem.}
  \end{center}
\end{figure}

In Figure \ref{fig:screwrb_err2} we compare the GRAC-2/3 method
(cf. \cite{2012-OrtZha-SINUM}) with the \BGFC scheme. We observe
numerical rates that are close to the predicted ones, and in
particular also the moderate improvement of \BGFC over GRAC-2/3
suggested in the previous section.

% \subsection{Anti-plane crack}
% %
% \label{sec:ext:crack}
% %

\subsection{Complex crystals}
\label{sec:ext:complex}
To formulate the \BGFC scheme for complex crystals, we return to the
point defect problem adressed in \S~\ref{sec:scheme}. 

\subsubsection{Atomistic model}
Each lattice site may now contain more than one atom (of the same or
different species). For simplicity suppose there are two atoms per
site, which we call species 1 and species 2. The deformation of the
lattice is now described by a deformation field $y : \L \to \R^d$ and
a shift $p : \L \to \R^d$. The deformed positions of species 1 are
given by $y_1(a) := y(a)$ and those of species 2 by $y_2(a) :=
y(a)+p(a)$, $a \in \L$. Let ${\bf y} := (y, p)$.  The site energy is
now a function
\begin{displaymath}
  \Phi_a({\bf y}) = \Phi_a\B( \b(y_i(b) - y_j(a)\b)_{\substack{i,j =
      1, 2 \\ b \in \L}} \B).
\end{displaymath}
% Note that this description, following \cite{2012-bvkort} is different
% from the more commonly employed kinematics in terms of a deformation
% $y = y_1$ and a shift $p = y_2-y_1$.

At present there exists no published regularity theory for complex
lattice defects corresponding to \cite{EhrOrtSha:defectsV1}, which we
employed in the discussion in \S~\ref{sec:scheme}, and the following
discussion is therefore based on unpublished notes \cite{unpub-bvkort,
  unpub-ort} and reasonable assumptions.

Let $\Phi_a^{\rm hom}$ be the site energy potential for the
defect-free lattice, then we assume that there is an equilibrium shift
$p_0$ such that ${\bf x} := (x, p_0)$ is a stable equilibrium
configuration. By this, we mean that, for all ${\bf v} = (v, r)$ with
$v, r : \Lhom \to \R^d$ compactly supported, 
\begin{align*}
  \sum_{a \in \L} \< \del \Phi_a({\bf x}), {\bf v} \> = 0 \quad
  \text{and} \quad 
  \sum_{a \in \L} \b\< \ddel \Phi_a({\bf x}) {\bf v}, {\bf v} \b\>
  \geq c_0 \b( |v|_{\Us}^2 + \| r \|_{\ell^2}^2 \b);
\end{align*}
that is, the configuration must also be stable under perturbations of
the shifts. This corresponds in fact to the classical notion of
stability in complex lattices; see \cite{E:2007a} and references
therein.

Then we define the energy difference functional
\begin{displaymath}
  \Ea({\bf u}) := \sum_{a \in \L} \Phi_a'({\bf u}), \qquad
  \text{where} \quad
  \Phi_a'({\bf u}) := \Phi_a({\bf x}+{\bf u}) - \Phi_a({\bf x}).
\end{displaymath}
It can again be shown that $\Ea$ is well-defined and regular on the
space \cite{unpub-ort}
\begin{displaymath}
{\bf U} := \b\{ {\bf v} = (v,r) : \L \to \R^{2d} \bsep v \in \Us,
  r \in \ell^2 \b\}.
\end{displaymath}
The exact atomistic problem now reads
\begin{displaymath}
  {\bf u}^\a \in \arg\min \b\{ \Ea({\bf u}) \bsep {\bf u} \in {\bf U} \b\}.
\end{displaymath}

\subsubsection{The \BGFC scheme} 
To define the BQCE and \BGFC schemes, we first define the Cauchy--Born
energy density by
\begin{displaymath}
  {\bf W}(\mF, p) := \det \mA^{-1} \Phi_0^{\rm hom}\b( (\mF x, p) \b),
\end{displaymath}
where $\Phi_a^{\rm hom}$ is the site energy potential for the
defect-free lattice. 
% \lz{LZ: shall we use $\Phi_0^{\rm hom}$ as the single crystal case?}
% CO: done
 Further, let $W'(\mG, q) := W(\mI+\mG, p_0+q) -
W(\mI, p_0)$.

Let the computational geometry be set up precisely as in
\S~\ref{sec:scheme:bqce}, and let 
\begin{displaymath}
  {\bf U}_h := \b\{ {\bf u}_h = (u_h, q_h) \bsep u_h, q_h \in \Ush \b\}.
\end{displaymath}
Note that, contrary to the usual practice, we require that both
displacement and shift are continuous functions. This is necessary in
order to be able to reconstruct atom positions. The BQCE energy
functional, as proposed in \cite{unpub-bvkort}, is given, for ${\bf
  u}_h = (u_h, q_h) \in {\bf U}_h$, by
\begin{displaymath}
  {\bf E}^{\rm b}({\bf u}_h) := \sum_{a \in \L} (1-\beta(a)) \Phi_a'({\bf u}_h)
  + \int_{\Omh} Q_h\b[ \beta W'(\D u_h, q_h)\b] \dx.
\end{displaymath}
Because of the loss of point symmetry in the interaction potential
there is also a reduction in the accuracy of the Cauchy--Born model
\cite{E:2007a} and in the blending scheme.  Indeed, the analysis in
\cite{unpub-bvkort} suggests that the best possible error that can be
expected for the complex lattice BQCE method is
\begin{displaymath}
  \| \D \overline{\ua} - \D u_h^{\rm b} \|_{L^2} + \| p^{\rm a} -
  p_h^{\rm b} \|_{\ell^2} \lesssim \| \D \beta \|_{L^2} + \text{best
    approx. err.} + \text{CB err.},
\end{displaymath}
that is the blended ghost force error now scales like $\| \D\beta
\|_{L^2}$. If $d = 2$, then it can be easily seen that $\| \D\beta
\|_{L^2} \gtrsim 1$ while for $d = 3$, one even gets $\| \D \beta
\|_{L^2} \gtrsim (\Ra)^{1/2}$. Thus, the standard BQCE scheme cannot
be optimised to become convergent in the energy-norm.

To formulate the \BGFC scheme, we renormalise $\Phi_a$ and $W$ a
second time,
\begin{align*}
  \Phi_a''({\bf u}) &:= \Phi_a({\bf x}+{\bf u}) - \Phi_a({\bf x}) - \<
  \del \Phi_a({\bf x}), {\bf u} \>, \quad \text{and} \\ 
  {\bf W}''(\mG, q) &:= W(\mI+\mG, p_0+q) - W(\mI, p_0) - \partial_\mF
  W(\mI, p_0) : \mG - \partial_{p} W(\mI, p_0) \cdot q,
\end{align*}
and define the \BGFC energy functional
\begin{displaymath}
  {\bf E}^{\rm rb}({\bf u}) := \sum_{a \in \L} (1-\beta(a)) \Phi_a''({\bf u})
  + \int_{\Omh} Q_h\b[ \beta W''(\D u, q)\b] \dx + \< {\bf L}^{\rm
    ren}, \bf u \>,
\end{displaymath}
where ${\bf L}^{\rm ren}$ is a suitable linear functional correcting
the forces in the defect core, defined analogously as $L^{\rm ren}$ in
\S~\ref{sec:scheme:rbqce}. The \BGFC scheme reads
\begin{equation}
  \label{eq:scheme:rbqce_complex}
  {\bf u}_h^{\rm rb} \in \arg\min \b\{ {\bf E}^{\rm rb}({\bf v}) \bsep
  {\bf v} \in {\bf U}_h \b\}.
\end{equation}

\subsubsection{Convergence rate}
While we leave a rigorous convergence theory for
\eqref{eq:scheme:rbqce_complex} to future work, we can still speculate
what rate of convergence may be expected.

Arguing analogously as in \S~\ref{sec:scheme:rbqce} and
\S~\ref{sec:analysis} we now observe that the error due to the blended
ghost forces can be bounded by
\begin{align*}
  \| \D \overline{\ua} - \D u_h^{\rm rb} \|_{L^2} + \| \overline{p^{\rm a}} -
  p_h^{\rm rb} \|_{L^2} &\lesssim \| \D \beta \|_{L^\infty} \b( \|\D
  u^\a \|_{L^2(\Omb)} + \| q^\a \|_{L^2(\Omb)} \b)   \\
  & \qquad + \text{best
    approx. err.} + \text{CB err.}.
\end{align*}

It is reasonable to expect that the regularity for the deformation
fields $y_1, y_2$ is similar as for simple lattices and therefore the
best approximation error is of the same order, i.e., $(\Ra)^{-d}$. The
Cauchy--Born modelling error can also be bounded by $\| \D^2
\widetilde{\ua} \|_{L^2(\R^d\setminus B_{\Ra})} + \| \D
\widetilde{p^\a} \|_{L^2(\R^d \setminus B_{\Ra})} \lesssim
(\Ra)^{-d}$, where we assumed again the same regularity for complex
lattice displacement fields as for the simple lattice case. (Since the
shift itself is a gradient, it is reasonable to assume that $\D
\widetilde{p^\a}$ decays like a second gradient.)

Thus, we are left to discuss the error due to the ghost
forces. Assuming the typical decay for point defects, $|\D
\widetilde{u^\a}| + |q^\a| \lesssim |x|^{-d}$ we obtain that, for a
quasi-optimal choice of $\beta$, satisfying $\|\D\beta\|_{L^\infty}
\lesssim (\Rb-\Ra)^{-1}$,
\begin{displaymath}
  \| \D \beta \|_{L^\infty} \b( \|\D
  u^\a \|_{L^2(\Omb)} + \| q^\a \|_{L^2(\Omb)} \b) \lesssim
  (\Rb-\Ra)^{-1} (\Ra)^{-d/2}.
\end{displaymath}
Upon choosing $\Rb-\Ra \approx \Ra$, this yields the rates
\begin{displaymath}
  \| \D \beta \|_{L^\infty} \b( \|\D
  u^\a \|_{L^2(\Omb)} + \| q^\a \|_{L^2(\Omb)} \b) \lesssim \cases{
    (\Ra)^{-2}, & d = 2, \\
    (\Ra)^{-5/2}, & d = 3,
  }
\end{displaymath}
which is the best approximation rate for $d = 2$ and a slightly
reduced rate for $d = 3$.

We therefore conclude that the expected rate of convergence for the complex lattice
\BGFC scheme, for point defects, is
\begin{displaymath}
  \| \D \overline{\ua} - \D u_h^{\rm rb} \|_{L^2} + \| \overline{p^{\rm a}} -
  p_h^{\rm rb} \|_{\ell^2} \lesssim \cases{ (\Ra)^{-2} \approx
    (\dof)^{-1}, & d = 2, \\
    (\Ra)^{-5/2} \approx (\dof)^{-5/6}, & d = 3.
  }
\end{displaymath}

With these heuristics in mind, we expect that it would be relatively
straightforward to generalize the analysis in \cite{LiOrShVK:2014} and
\S~\ref{sec:analysis} and thus obtain the first rigorously convergent
a/c coupling scheme for complex crystals.

\section{Analysis}
\label{sec:analysis}
\def\Pdef{\mathscr{P}^{\rm def}}
%
% \subsection{The simple lattice case}
% \label{sec:analysis_simple}
%
For our rigorous error estimates we focus on a simplified point defect
problem, following \cite{LiOrShVK:2014}. We shall cite several results
that are summarized in \cite{LiOrShVK:2014} but drawn from other
sources, but for the sake of convenience we will only cite
\cite{LiOrShVK:2014} as a reference. A range of generalizations are
possible but require some additional work, and in particular a more
complex notation.

We assume $\L \equiv \Lhom$, with globally homogeneous site energies
with a finite interaction radius in {\em reference
  configuration}. That is, we assume that there exists $\Rg \subset
B_{\rcut} \cap (\L \setminus \{0\})$, finite, and $V \in
C^4((\R^d)^\Rg)$ such that $\Phi_a(y) = V(Dy(a))$, where $Dy(a) :=
(D_\rho y(a))_{\rho \in \Rg}$ and $D_\rho y(a) := y(a+\rho) -
y(a)$. We assume throughout that $V$ is point symmetric, i.e.,
$-\Rg=\Rg$ and $V( (-g_{-\rho}) ) = V((g_\rho))$ for all
$(g_\rho)_{\rho \in \Rg} \in (\R^d)^\Rg$.

A defect is introduced by adding an external potential $\Pdef \in
C^2(\Us)$, which only depends on $Du(a), |a| < \Rdef$. The atomistic
problem now reads
\begin{equation}
  \label{eq:ana:atm}
  \ua \in \arg\min \b\{ \Ea(v) + \Pdef(v) \bsep v \in \Us \b\}.
\end{equation}
We call a point $\ua$ a strongly stable solution to \eqref{eq:ana:atm}
if there exists $\gamma > 0$ such that
\begin{displaymath}
  \< \del[\Ea + \Pdef](\ua), v \> = 0 \quad \text{and} \quad
  \< \ddel[\Ea+\Pdef](\ua) v, v \> \geq \gamma |v|_{\Us}^2 
  \quad \forall v \in \Us.
\end{displaymath}

The \BGFC approximation to \eqref{eq:ana:atm}  is given by
\begin{equation}
  \label{eq:ana:rbqce}
  \ubr \in \arg\min \b\{ \Ebr(v_h) + \Pdef(v_h) \bsep v_h \in \Ush
  \b\},
\end{equation}
using the notations of \S~\ref{sec:scheme:rbqce}. 

\subsection{Additional assumptions}
\label{sec:simple:assumptions}
\def\TT{\mathcal{T}}
\def\Oma{\Omega^{\rm a}}
\def\Omb{\Omega^{\rm b}}
\def\Omc{\Omega^{\rm c}}
We now summarize the main assumptions we require to state our rigorous
convergence results. All assumptions can be satisfied in practice, and
are discussed in detail in \cite{LiOrShVK:2014}.

We assume that $\beta \in C^{2,1}(\R^d)$, $0 \leq \beta \leq 1$. Let
$\rcut' := 2 \rcut + \sqrt{d}$, $\Oma := {\rm supp}(1-\beta) +
B_{\rcut'}$, $\Omb := {\rm supp}(\D\beta) + B_{\rcut'}$ and $\Omc :=
({\rm supp}(\beta) + B_{\rcut'}) \cap \Omh$. Then we require that
there exist radii $\Ra \leq \Rb \leq \Rc$ and constants $C_{\rm b},
C_\Omega$ such that
\begin{align*}
  & \Rb \leq C_{\rm b} \Ra, \qquad \| \D^j \beta \|_{L^\infty} \leq
  C_{\rm b} (\Ra)^{-j}, j = 1, 2, 3;  \\
  & {\rm supp}(\beta) \supset
  B_{\Ra+\rcut'}, \qquad {\rm supp}(1-\beta) \subset B_{\Rb-\rcut'};
  \\
  &B_{\Rc} \subset \Oma\cup\Omc \qquad \text{and} \qquad \Rc \geq
  C_{\Omega} (\Ra)^2.
\end{align*}

To state the final assumption that we require on $\Th$, we first need
to define a piecewise affine interpolant of lattice functions.

The first is a piecewise affine interpolant. If $d = 2$, let
$\hat{\mathcal{T}} := \{ \hat{T}_1, \hat{T}_2 \}$, where $\hat{T}_1 =
{\rm conv}\{0, e_1, e_2\}$ and $\hat{T}_2 = {\rm conv}\{e_1, e_2,
e_1+e_2\}$, where ${\rm conv}$ denotes the convex hull of a set of
points.
% \lz{LZ: {\rm conv} not defined yet}.
If $d = 3$, let $\hat{\mathcal{T}} := \{ \hat{T}_1, \dots, \hat{T}_6
\}$ be the standard subdivision of $[0, 1]^3$ into six tetrahedra (see
\cite[Fig. 1]{LiOrShVK:2014}). Then $\TT := \bigcup_{\ell \in \L}
(\ell+\hat{\mathcal{T}})$ defines a regular and uniform triangulation
of $\R^d$ with node set $\L$. For each $v : \L \to \R^m$, there exists
a unique $\barv \in C(\R^d; \R^m)$ such that $\barv(\ell) = v(\ell)$
for all $\ell \in \L$. In particular, we note that the semi-norms $\|
\D \barv \|_{L^2}$ and $|v|_{\Us}$ are equivalent.

Our final requirement on the approximation parameters is that there
exists a constant $C_h$ such that
\begin{align*}
  &\Th \cap \Oma \equiv \TT \cap \Oma \qquad \qquad
  \max_{T \in \Th} h_T^d / |T| \leq C_h, \\
  &h(x) \leq C_h \max\b(1, |x|/\Ra\b) \qquad \text{and} \qquad \#\Th
  \leq C_h (\Ra)^2 \log(\Ra).
\end{align*}
By $\Th \cap \Oma \equiv \TT \cap \Oma$ we mean that, if $T \in \TT,
T\cap\Oma \neq \emptyset$ then $T \in \Th$ (and hence also
vice-versa). The condition $h(x) \leq C_h \max\b(1, |x|/\Ra\b)$ can be
weakened; see e.g. \cite{2012-CMAME-optbqce, 2012-MATHCOMP-qce.pair}.

We remark that $(\beta, \Th)$ are the main approximation parameters,
while the ``regularity constants'' ${\bf C} = (C_{\rm b}, C_h,
C_\Omega)$ are ``derived parameters''. In the following we {\em fix}
the constants ${\bf C}$ to some given bounds, and admit any pair
$(\beta, \Th)$ satisfying the foregoing conditions with these
constants. When we write $A \lesssim B$, then we mean that there
exists a constant $C$ depending only on ${\bf C}$ (as well as on the
solution and on the model) but not on $(\beta, \Th)$ such that $A \leq
C B$.

\subsection{Convergence result}
\label{sec:simple:conv}
The following convergence result is a direct extension of Theorem 3.1,
Proposition 3.2 and Theorem 3.3 in \cite{LiOrShVK:2014} to the \BGFC
method. The proof of the theorem is given in the next two sections.

\begin{theorem}
  \label{th:simple}
  Let $\ua$ be a strongly stable solution to \eqref{eq:ana:atm}, then
  for any given set of constants ${\bf C}$ there exist $C, C', R^\a_0 > 0$
  such that, for all $(\beta, \Th)$ satisfying the conditions of
  \S~\ref{sec:simple:assumptions}, and in addition $\Ra \geq \Ra_0$,
  there exists a solution $\ubr$ to \eqref{eq:ana:rbqce} such that
  \begin{align}
    \label{eq:simple:normest}
    \| \D \overline{\ua} - \D \ubr \|_{L^2} &\leq C (\Ra)^{-d/2-1} \leq
    C' \b(\smfrac{\log \#\Th}{\#\Th} \b)^{1/2+1/d}, 
    \qquad \text{and} \\
    \label{eq:simple:energyest}
    \b| \Ea(\ua) - \Ebr(\ubr) \b| &\leq C (\Ra)^{-d-2} 
    \leq C'  \b(\smfrac{\log \#\Th}{\#\Th} \b)^{1+2/d}.
  \end{align}
\end{theorem}

\begin{remark}
  1. The only assumption we made that represents a genuine restriction
  of generality is that $\|\D^j \beta\|_{L^\infty} \lesssim
  (\Ra)^{-j}$. We require this to prove stability of the BQCE and BGFC
  schemes.

  2. However, the proof of \eqref{eq:simple:energyest} shows that the
  energy error would be sub-optimal if we chose a narrower blending
  region (and thus a slower rate of $\|\D^j \beta\|_{L^\infty}
  \lesssim (\Ra)^{-sj}$ for some $s < 1$). This appears to contradict
  our numerical results in \S~\ref{sec:numerical_tests} and suggests
  that the energy error estimate may be suboptimal.
\end{remark}

\subsection{Proof of the energy norm error estimate}
\label{sec:simple:prf}
To prove the result we will need to refer to another technical tool
from \cite{LiOrShVK:2014}, namely a $C^{2,1}$-conforming
multi-quintic, which we use to measure the regularity of an atomistic
displacement.

For $v : \L \to \R^m$ and $i = 1, \dots, d$, let ${\sf d}_i^0 v(\ell)
:= v(\ell)$; ${\sf d}_i^1 v(\ell) := \frac12 (u(\ell+e_i) -
u(\ell-e_i))$ and ${\sf d}_i^2 v(\ell) := u(\ell+e_i) - 2 u(\ell) +
u(\ell-e_i)$. Lemma 2.1 in \cite{LiOrShVK:2014} states that, for each
$\ell \in \L$ there exists a unique multi-quintic function $\tilv :
\ell+[0,1]^d \to \R^m$ defined through the conditions
\begin{displaymath}
  \partial_{x_1}^{\alpha_1} \cdots \partial_{x_d}^{\alpha_d}
  \tilv(\ell') = {\sf d}_1^{\alpha_1} \cdots {\sf d}_{d}^{\alpha_d}
  v(\ell') \qquad \forall \ell' \in \ell+\{0,1\}^d,  \alpha \in
  \{0,1,2\}^d, \|\alpha \|_\infty \leq 2,
\end{displaymath}
and moreover, that the resulting piecewise defined function on $\R^d$
belongs to $\tilv \in C^{2,1}(\R^d; \R^m)$.

We begin the proof of Theorem \ref{th:simple} by noting that the
renormalised site energy potential
\begin{displaymath}
  V''(Du) := V(\Rg+Du)-V(\Rg)-\< \del V(\Rg), Du \>
\end{displaymath}
is an admissible potential for
\cite[Thm. 3.1]{LiOrShVK:2014}. Further, the conditions we put forward
in \S~\ref{sec:simple:assumptions} are precisely those we need to
apply \cite[Thm. 3.1]{LiOrShVK:2014} with $V$ replaced with $V''$,
thus treating \BGFC as a simple BQCE method. Hence, we obtain that,
under the conditions of Theorem \ref{th:simple}, there exists a
solution $\ubr$ to \eqref{eq:ana:rbqce} and constants $C_1, C_2$
depending only on ${\bf C}$, but independent of the approximation
parameters, such that
\begin{align}
  \notag
  \| \D \overline{\ua} - \D \ubr \|_{L^2} \leq C_1 \| \D^2 \beta
  \|_{L^2} + C_2 \B(&~ \| \D \bar{u}^\a \|_{L^2(\R^2\setminus
    B_{\Rc/2})} + \| h \D^2 \tilde{u}^\a \|_{L^2(\Omc)} \\
  \label{eq:simpleprf:errest1}
  & + \| h^2 \D^3 \tilde{u}^\a
  \|_{L^2(\Omc)} + \| \D^2 \tilde{u}^\a \|_{L^4(\Omc)}^2 \B).
\end{align}
Here we did not write out some terms that are trivially dominated by
those that we did write. In the following we write $u \equiv u^\a$.

The group preceded by the constant $C_2$ cannot be further improved,
but we will analyse in more detail the group $C_1
\|\D^2\beta\|_{L^2}$. This term arises from the coarsening and
modelling error analysis of the BQCE scheme in \S~6.1 and \S~6.2 of
\cite{LiOrShVK:2014} (see also the summary in \S 4.3 of
\cite{LiOrShVK:2014}). We can replace $C_1 \|\D^2\beta\|$ in
\eqref{eq:simpleprf:errest1} with two terms $\epsilon_\beta^{\rm
  coarse} + \epsilon_\beta^{\rm model}$, which we discuss next, to
obtain
\begin{align}
  \notag
  \| \D \baru - \D \ubr \|_{L^2} \leq \epsilon_\beta^{\rm
  coarse} + \epsilon_\beta^{\rm model} + C_2 \B(&~ \| \D \baru \|_{L^2(\R^2\setminus
    B_{\Rc/2})} + \| h \D^2 \tilu \|_{L^2(\Omc)} \\
  \label{eq:simpleprf:errest2}
  & + \| h^2 \D^3 \tilu
  \|_{L^2(\Omc)} + \| \D^2 \tilu \|_{L^4(\Omc)}^2 \B).
\end{align}

From Lemma 6.1 and Lemma 6.2 of \cite{LiOrShVK:2014} it can be readily
seen that the coarsening error contribution to $C_1
\|\D^2\beta\|_{L^2}$ is 
\begin{displaymath}
  \epsilon_\beta^{\rm coarse} \lesssim \| \partial W''(\mI+ \D\tilu) \D^2
  \beta \|_{L^2}.
\end{displaymath}
Using that fact that $\partial W''(\mI)  = 0$ (due to the
renormalisation, the reference configuration is now stress free) we
therefore obtain
\begin{equation}
  \label{eq:simpleprf:epscoarse}
  \epsilon_\beta^{\rm coarse} \lesssim \b\| \b[ \partial W''(\mI+
  \D\tilu) - \partial W''(\mI) \b] \D^2
  \beta \b\|_{L^2} 
  \lesssim \b\| \,|\D \tilu|\,|\D^2\beta|\,\b\|_{L^2}.
\end{equation}

To obtain a more precise control on the modelling error contribution
we first take a closer look at the term ${\rm T}_1$ defined in
Equation (6.6) in \cite{LiOrShVK:2014},
\begin{displaymath}
  {\rm T}_1 = \beta \partial W'' - \sum_{\ell \in \L}
  \beta(\ell) \sum_{\rho \in \Rg} \b[ V_{,\rho}'' \otimes \rho\b] \omega_\rho(\ell-x),
\end{displaymath}
where $\partial W'' = \partial W''(\mI+\D\tilu(x))$, $V_{,\rho}'' =
V_{,\rho}''((\mI+\D\tilu(x))\Rg)$ and $\omega_\rho(y) = \int_{s = 0}^1
\bar{\zeta}(y+s\rho) \ds$, with $\bar{\zeta}$ being the P1 hat
function for the origin on the mesh $\TT$. The $O(\|\D^2\beta\|)$ term
now arises by expanding $\beta$,
\begin{displaymath}
  \beta(\ell) = \beta(x) + \D\beta(x) \cdot (\ell-x) + R_\beta(x;\ell),
\end{displaymath}
with $|R_\beta(x;\ell)| \lesssim \|\D^2\beta\|_{L^\infty}$ due to the
fact that only $\ell$ within a fixed radius around $x$ are considered
in the seemingly infinite sum. Thus, exploiting also the fact that
$V''_{,\rho}(\Rg) \equiv 0$ and $\sum_{\ell \in \L}
\omega_\rho(\ell-x) = 1$ (see \cite[Eq. (4.18)]{LiOrShVK:2014}), and
denoting the terms that lead to the second group in
\eqref{eq:simpleprf:errest2} by ``$\cdots$'', we obtain
\begin{align*}
  |{\rm T}_1| &\leq \|\D^2\beta\|_{L^\infty} \sum_{\rho\in\Rg} \B|
  \b(V_{,\rho}''(\Rg+\D\tilu) - V_{,\rho}''(\Rg)\b) \otimes \rho \sum_{\ell
    \in \L} \omega_\rho(\ell-x) \B| + \dots \\
  &\lesssim \|\D^2\beta\|_{L^\infty} \sum_{\rho\in\Rg} \b|
  V_{,\rho}''(\Rg+\D\tilu) - V_{,\rho}''(\Rg) \b| + \dots \\
  &\lesssim \|\D^2 \beta \|_{L^\infty} |\D\tilu(x)| + \dots.
\end{align*}
The improved modelling error estimate is obtained by taking the
$L^2$-norm of ${\rm T}_1$, interpreting it as a function of $x$. Upon
noting that ${\rm T}_1(x) = 0$ outside $\Omb$, we therefore obtain
\begin{equation}
  \label{eq:simpleprf:epsmodel}
  \epsilon_\beta^{\rm model} \lesssim \| \D^2\beta\|_{L^\infty} \|
  \D\tilu \|_{L^2(\Omb)}.
\end{equation}

Combining \eqref{eq:simpleprf:errest2}, \eqref{eq:simpleprf:epscoarse}
and \eqref{eq:simpleprf:epsmodel} we arrive at
\begin{align}
  \notag
  \| \D \baru - \D \ubr \|_{L^2} &\leq C \B(\|\D^2\beta\|_{L^\infty}
  \|\D\tilu \|_{L^2(\Omb)} + \| \D \baru \|_{L^2(\R^2\setminus
    B_{\Rc/2})} + \| h \D^2 \tilu \|_{L^2(\Omc)} \\
  \label{eq:simpleprf:errest3}
  & \qquad + \| h^2 \D^3 \tilu
  \|_{L^2(\Omc)} + \| \D^2 \tilu \|_{L^4(\Omc)}^2 \B),
\end{align}
for a constant $C$ that depends only on ${\bf C}$, but is independent
of $(\beta, \Th)$.

Following the proof of \cite[Thm. 3.3]{LiOrShVK:2014} in \S 3.2.2 of
\cite{LiOrShVK:2014} it is straightforward now to obtain the rate
\eqref{eq:simple:normest}. Towards its proof we only remark that,
according to \cite[Lemma 2.3]{LiOrShVK:2014}, $|\D\tilu(x)|\lesssim
|x|^{-d}$ and hence, using the assumption $\|\D^2\beta\|_{L^\infty}
\lesssim (\Ra)^{-2}$, we obtain
\begin{displaymath}
  \|\D^2\beta\|_{L^\infty} \, \|\D\tilu \|_{L^2(\Omb)} \lesssim
  (\Ra)^{-2} \bg( \int_{\Ra}^{\Rb} r^{d-1} r^{-2d} \dr \bg)^{1/2}
  \lesssim (\Ra)^{-2-d/2}.
\end{displaymath}

This completes the proof of \eqref{eq:simple:normest}.

\subsection{Proof of the energy error estimate}
As in the case of the energy-norm error estimate we only modify some
specific parts of the proof for the BQCE case in \S~6.3 of
\cite{LiOrShVK:2014}, as required to obtain the improved energy error
estimate. To follow the notation let $u \equiv \ua$ and $u_h \equiv
\ubr$. Further, we recall that $\Pi_h u \in \Ush$ is a
best-approximation of $u$. For the following proof we do not need to
know its precise definition, but only remark that $\Pi_h u = u$ in
$\Oma$ and, as an intermediate step in the proof of
\eqref{eq:simple:normest} one obtains
\begin{displaymath}
  \| \D \tilu - \D \Pi_h u \|_{L^2(\Omc)} + \| \D u_h - \D \Pi_h u
  \|_{L^2} \lesssim (\Ra)^{-1-d/2}.
\end{displaymath}

\def\tilE{\tilde{\mathcal{E}}}

Following the proof of the BQCE energy error estimate in \S~6.3 of
\cite{LiOrShVK:2014}, we split the energy error into $\Ea(u) -
\Ebr(u_h) = T_1 + T_2 + T_3$, where 
\begin{align*}
  &T_1 = \Ea(u) - \tilE, \qquad T_2 = \tilE - \Ebr(\Pi_h u), \qquad
  T_3 = \Ebr(\Pi_h u) - \Ebr(u_h), \\
  &\text{and } \tilE = \sum_{\ell \in \L} (1-\beta(\ell)) V''(Du(\ell)) +
  \int_{\R^d} [Q_h\beta] W''(\D \tilu) \dx.
\end{align*}
We treat the terms in the same order as in \cite{LiOrShVK:2014}.

As in \cite{LiOrShVK:2014}, the term $T_3$ can be bounded by
\begin{equation}
  \label{eq:enprf:T3}
  |T_3| \lesssim \| \D u_h - \D \Pi_h u
  \|_{L^2} \lesssim (\Ra)^{-2-d}.
\end{equation}

The term $T_1$ is split 
further into
\begin{align*}
  T_1 &= \sum_{\ell \in \L} \beta(\ell) \b( V''(Du(\ell)) - W''(\D
  \tilu(\ell)) \b) \\
  & \qquad + \int_{\R^d} \b( [Q_h\beta] W''(\D\tilu) -
  I_1[\beta W''(\D\tilu)] \b) \dx 
  =: T_{1,1} + T_{1,2},
\end{align*}
where $I_1$ denotes the P1 nodal interpolant for the atomistic mesh
$\TT$. The second term is essentially a quadrature error and following
the proof of (6.12) in \cite{LiOrShVK:2014} (but note that the
inverse-estimate trick is not required in our present setting) it is
easy to see that
\begin{displaymath}
  |T_{1,2}| \lesssim \sum_{T \in \TT} \| \D^2 [\beta W''(\D \tilu)] \|_{L^\infty(T)}.
\end{displaymath}
The summand vanishes, unless $T \subset \Omc$. In the latter case, we
have
\begin{align*}
  \| \D^2 [\beta W''(\D \tilu)] \|_{L^\infty(T)} &\lesssim \|
  \D^2\beta \|_{L^\infty} \| \D\tilu \|_{L^2(T)}^2 + \|
  \D\beta\|_{L^\infty(T)} \| \partial W''(\D\tilu) \D^2 \tilu
  \|_{L^\infty(T)}  \\
  &\qquad + \| \partial W''(\D\tilu) \D^3 \tilu
  \|_{L^\infty} + \| \partial^2 W''(\D\tilu) \D^2 \tilu
  \|_{L^\infty(T)}^2 \\
  &\lesssim \|
  \D^2\beta \|_{L^\infty} \| \D\tilu \|_{L^2(T)}^2 + \|
  \D\beta\|_{L^\infty} \| \D\tilu \|_{L^2(T)} \|\D^2\tilu
  \|_{L^2(T)} \\
  & \qquad + \| \D \tilu \|_{L^2(T)} \| \D^3 \tilu \|_{L^2(T)} + \|
  \D^2\tilu \|_{L^2(T)}^2.
\end{align*}
Therefore, we obtain
\begin{align}
  \notag
  |T_{1,2}| &\lesssim \|
  \D^2\beta \|_{L^\infty} \| \D\tilu \|_{L^2(\Omb)}^2 + \|
  \D\beta\|_{L^\infty} \| \D\tilu \|_{L^2(\Omb)} \|\D^2\tilu
  \|_{L^2(\Omb)} \\
  \label{eq:enprf:T12}
  & \qquad + \| \D \tilu \|_{L^2(\Omc)} \| \D^3 \tilu \|_{L^2(\Omc)} + \|
  \D^2\tilu \|_{L^2(\Omc)}^2 \\
  \notag
  &\lesssim (\Ra)^{-2-d} + (\Ra)^{-2-d} + (\Ra)^{-2-d} +
  (\Ra)^{-2-d} \approx (\Ra)^{-2-d}.
\end{align}

To estimate $T_{1,1}$ we begin by noting that
\begin{align*}
  V''(Du) - W''(\D u) &= V''(\Rg+Du) - V''(\Rg+\D_\Rg u) \\
  &= \< \del V''(\Rg+\D_\Rg u), Du - \D_\Rg u \> \\
  & \qquad 
  + \< \ddel V(\Theta) (Du-\D_\Rg u), Du - \D_\Rg u \>,
\end{align*}
where $\Theta \in {\rm conv}\{\Rg+Du, \Rg+\D_\Rg u\}$ and $\D_\Rg
  u := (\D_\rho u)_{\rho \in \Rg}$.
% \lz{LZ: $\D_\Rg u$ ($\D u \cdot\Rg$?) not defined}.
Then, continuing to argue as in Lemma 6.5 and \cite{LiOrShVK:2014}
(performing a Taylor expansion on $Du$ and exploiting the point
symmetry of $V''$, and in particular exploiting the fact that $\|\del
V''(\Rg+\D_\Rg u)\| \lesssim |\D\tilu|$) we obtain
\begin{align*}
  \b| \< \del V''(\Rg+\D_\Rg u), Du - \D_\Rg u \> \b| &\lesssim | |\D
  u| \| \D^3 \tilu \|_{L^\infty(\nu_x)}, \quad \text{and} \\
\b|\< \ddel V(\Theta) (Du-\D_\Rg u), Du - \D_\Rg u \>\b| &\lesssim \|
\D^2 \tilu \|_{L^\infty(\nu_x)}^2,
\end{align*}
where $\nu_x = B_{\rcut'}(x)$. Using the fact that $\D^3\tilu$ is
a piecewise polynomial we can use the inverse inequalities in (5.7) of
\cite{LiOrShVK:2014} to obtain, as in (6.14) of \cite{LiOrShVK:2014},
that
\begin{align}
  \notag
  |T_{1,1}| &\lesssim \| \D \tilu \|_{L^2(\Omc)} \| \D^3\tilu
  \|_{L^2(\Omc)} + \| \D^2 \tilu \|_{L^2(\Omc)}^2  \\
  \label{eq:enprf:T11}
  &\lesssim (\Ra)^{-2-d} + (\Ra)^{-2-d} \approx (\Ra)^{-2-d}.
\end{align}

Since $\Pi_h u = u$ in $\Oma$, the term $T_2$ simplifies to
\begin{displaymath}
  T_2 = \int_{\R^d} [Q_h\beta] \b( W''(\D\tilu) - W''(\D\tilu
  - \D e) \b) \B) \dx,
\end{displaymath}
where $e := \tilu - \Pi_h u$ and we used the fact that $Q_h[\beta
W''(\D\Pi_h u)] = [Q_h\beta] W''(\D\Pi_h u)$. We begin by expanding $W''(\D\tilu) - W''(\D\tilu
  - \D e) = \partial W''(\D\tilu) : \D e + O(|\D e|^2)$ to obtain
\begin{align*}
  T_2 &\lesssim \bg| \int_{\R^d} [Q_h\beta] \b[ \partial W''(\D\tilu)
  : \D e \b] \dx\bg| + \| \D e \|_{L^2(\Omc)}^2 \\
  &\lesssim \bg| \int_{\R^d} \beta \b[ \partial W''(\D\tilu)
  : \D e \b] \dx\bg| +
  \bg| \int_{\R^d} [Q_h\beta - \beta] \b[ \partial W''(\D\tilu)
  : \D e \b] \dx\bg|  \\
  & \qquad +
  \| \D e \|_{L^2(\Omc)}^2 \lesssim T_{2,1}+T_{2,2} + T_{2,3}.
\end{align*}
To treat $T_{2,1}$ we integrate by parts, and then use $|\partial
W''(\D\tilu)| = |\partial W(\mI+\D\tilu) - \partial W(\mI)| \lesssim
|\D\tilu|$, $\|e \|_{L^2(T)} \lesssim \| h^2 \D^2 \tilu \|_{L^2(T)}$,
and $\|h^2 \D^2\tilu \|_{L^2(\Omc)} \lesssim (\Ra)^{-1-d/2}$ to
estimate
\begin{align*}
  T_{2,1} &= \bg| \int_{\R^d} - {\rm div} \B( \beta \partial
  W''(\D\tilu) \B) \cdot e \dx\bg| \\
  &\lesssim \| \D\beta \|_{L^\infty} \| \partial
  W''(\D\tilu) \|_{L^2(\Omc)} \| e \|_{L^2(\Omc)} + \| \D^2\tilu
  \|_{L^2(\Omc)} \| e \|_{L^2(\Omc)} \\
  &\lesssim \| \D\beta \|_{L^\infty} \| \D\tilu \|_{L^2(\Omc)} \| e \|_{L^2(\Omc)} + \| \D^2\tilu
  \|_{L^2(\Omc)} \| e \|_{L^2(\Omc)} \\
  &\lesssim (\Ra)^{-1} (\Ra)^{-d/2} (\Ra)^{-1-d/2} + (\Ra)^{-1-d/2}
  (\Ra)^{-1-d/2} \approx (\Ra)^{-2-d}.
\end{align*}
The terms $T_{2,2}$ and $T_{2,3}$ are estimated analogously, by
$T_{2,j} \lesssim (\Ra)^{-2-d}$, $j = 2, 3$, and thus we obtain that
\begin{equation}
  \label{eq:enprf:T2}
  |T_2| \lesssim (\Ra)^{-2-d}.
\end{equation}

% Estimating
% % \begin{align*}
% %   \b|W''(\D \tilu) - W''(\D\tilu - \D e)\b| 
% %   &\lesssim \b|[\partial W''(\D\tilu) - \partial W''(\mI)] : \D e \b|
% %   + |\D e|^2  \\
% %   & \lesssim |\D\tilu| |\D e| + |\D e|^2,
% % \end{align*}
% % we readily obtain
% % \begin{align}
% %   \notag 
% %   |T_{2,1}| &\lesssim \| \D \tilu \|_{L^2(\Omc)} \| \D e \|_{L^2(\Omc)} + \| \D
% %   e \|_{L^2(\Omc)}^2 \\
% %   \label{eq:enprf:T2}
% %   &\lesssim (\Ra)^{-1-d} + (\Ra)^{-2-d} 
% %   \lesssim (\Ra)^{-1-d}.
% % \end{align}
% % ***!***

% \begin{align*}
%   \b|W''(\D \tilu) - W''(\D\tilu - \D e)\b| 
%   &\lesssim \b|\partial W''(\D\tilu) : \D e \b|
%   + |\D e|^2  \\
%   & \lesssim \b| \b[\partial W(\mI+\D\tilu) - \partial W(\mI)\b] : \D e \b| + |\D e|^2,
% \end{align*}
% we readily obtain
% \begin{align}
%   \notag 
%   |T_{2,1}| &\lesssim \| \D \tilu \|_{L^2(\Omc)} \| \D e \|_{L^2(\Omc)} + \| \D
%   e \|_{L^2(\Omc)}^2 \\
%   \label{eq:enprf:T2}
%   &\lesssim (\Ra)^{-1-d} + (\Ra)^{-2-d} 
%   \lesssim (\Ra)^{-1-d}.
% \end{align}
% ***!***

Combining \eqref{eq:enprf:T3}, \eqref{eq:enprf:T12},
\eqref{eq:enprf:T11} and \eqref{eq:enprf:T2} completes the proof of
\eqref{eq:simple:energyest}.

\bibliographystyle{plain}
\bibliography{qc}

\end{document}

%% file: notation.tex
\renewcommand{\cases}[1]{\left\{ \begin{array}{rl} #1 \end{array} \right.}
\newcommand{\smfrac}[2]{{\textstyle \frac{#1}{#2}}}

\def\XXint#1#2#3{{\setbox0=\hbox{$#1{#2#3}{\int}$ }
\vcenter{\hbox{$#2#3$ }}\kern-.6\wd0}}
\def\b{\big}
\def\B{\Big}
\def\bg{\bigg}

\def\sep{\,|\,}
\def\bsep{\,\b|\,}

\def\defc{\rm def}

\def\R{\mathbb{R}}

\def\Z{\mathbb{Z}}

\def\dx{\,{\rm d}x}

\def\dr{\,{\rm d}r}

\def\ds{\,{\rm d}s}

\def\pp{\partial}

\def\<{\langle}
\def\>{\rangle}

\def\mA{{\sf A}}
\def\mB{{\sf B}}
\def\mF{{\sf F}}
\def\mG{{\sf G}}

\def\mI{{\sf I}}

\def\mQ{{\sf Q}}

\def\D{\nabla}
\def\del{\delta}
\def\ddel{\delta^2}

\renewcommand\a{{\rm a}}
\renewcommand\c{{\rm c}}

\def\L{\Lambda}
\def\Lhom{\L^{\rm hom}}

\newcommand{\Omb}{{\Omega^{\beta}}}

\def\Us{\mathscr{U}}

\def\E{\mathscr{E}}
\def\Ea{\E^\a}

\def\Omh{{\Omega_h}}
\def\Omc{\Omega^\c}

\def\tilu{\tilde u}
\def\tilv{\tilde v}

\def\barv{\bar v}

\def\baru{\bar u}

\def\TT{\mathscr{T}}

\def\Th{\mathscr{T}_h}
\def\Nh{\mathscr{N}_h}